\documentclass [a4paper,twoside,12pt]{article}

\usepackage{amsmath,amsfonts,amssymb}
\usepackage{vmargin,graphicx,theorem}
\usepackage{mathtools}
\usepackage{enumerate}
\usepackage[english]{babel}
\usepackage{enumitem}
\setpapersize[portrait]{A4}
\usepackage{hyperref}
\setmarginsrb{1cm}{2.2cm}{1cm}{1.7cm}{0cm}{0.1cm}{0.5cm}{2cm}

\newcommand{\R}{\ensuremath{\mathbb{R}}}

\newtheorem{ethm}{}
\newtheorem{theo}[ethm]{Theorem}

\newtheorem{ecor}[ethm]{Corollary}

\newtheorem{eprop}[ethm]{Proposition}

\newtheorem{elem}[ethm]{Lemma}

\newtheorem{erem}[ethm]{Remark}

\newcommand{\proofend}{~\hfill$\square$}

\newenvironment{eproof}[1][]
{%
  \noindent
  \emph{\textbf{Proof%
  \if\relax\detokenize{#1}\relax
  \else\ of #1%
  \fi}}\\
}
{%
  \proofend\\
}

\newcommand{\abs}[1]{\ensuremath{{\left| #1 \right|}}} 
\newcommand{\NRM}[1]{\ensuremath{{\left\Vert #1\right\Vert}}} 
\newcommand{\CLSI}{\ensuremath{\mathrm{C}_{\mathrm{LSI}}}}
\newcommand{\Const}{\ensuremath{\mathrm{C}}}
\newcommand{\Constinit}{\ensuremath{\mathrm{C_{\mathrm{init}}}}}
\newcommand{\Hyp}{\ensuremath{(\mathrm{H}')}}
\newcommand{\Hypp}{\ensuremath{(\mathrm{H})}}

\renewcommand{\phi}{\varphi}

\newcommand{\Ent}[2]{\mathbf{Ent}( #1 \,,\, #2 )}

\newcommand{\diver}[1]{\nabla \cdot #1}

\newcommand{\grad}{\nabla}
\newcommand{\E}[1]{\mathbb{E}\left[#1\right]}
\newcommand{\N}{\mathbb{N}}
\newcommand{\lap}{\Delta}

\begin{document}

\title{\sl Entropic Propagation and Generation of Chaos for McKean-Vlasov Diffusions with Polynomial Growth}
\author{ Matteo Miannay}
\date{July 10, 2026}
\maketitle 
\begin{center}
\small{ Universit\'e de Rennes, ENS de Rennes,
IRMAR, UMR CNRS 6625\\
matteo.miannay@ens-rennes.fr}
\end{center}

\bigskip

\begin{center}
{\bf Abstract}
\end{center}

We study time-uniform \textit{propagation} and \textit{generation} of chaos for a class of McKean-Vlasov equations with polynomial growth potentials. Using an entropy dissipation method based on time-uniform logarithmic Sobolev inequalities, we derive uniform-in-time bounds depending on the regularity of the potentials and the assumptions on the initial conditions. We thereby not only establish a result stated but unproven in \cite{malrieu}, but also extend it to non-chaotic and non-exchangeable initial conditions, moving beyond the classical Lipschitz setting at the expense of a logarithmic loss in the convergence rate.\selectlanguage{english}

\bigskip
\section*{Framework and overview of the results}



Given smooth potentials $V,W : \R^d \longrightarrow \R$ and a family  $(B^i)_{1 \leq i \leq N}$ of independent Brownian motions on $\R^d$, we consider the following particle system in $\R^{dN}$:
\begin{equation}\label{systemedeparticules}
\begin{cases}
\displaystyle    \mathrm{d}X_t^{i,N} = \sqrt{2}\,\mathrm{d}B_t^i \;-\; \nabla V\bigl(X_t^{i,N}\bigr)\,\mathrm{d}t  \;-\; \frac{1}{N}\sum_{j=1}^N \nabla W\bigl(X_t^{i,N} - X_t^{j,N}\bigr)\,\mathrm{d}t, \quad 1 \leq i \leq N \\
    \mathcal{L}(X_0^{1,N},\cdots,X_0^{N,N}) = u_0^{(N)}.
\end{cases}
\end{equation}
We define $u_t^{(N)}$ as the law at time $t$ of the solution of \eqref{systemedeparticules} in $\R^{dN}$.
In terms of particle dynamics, the potential $V$ is a confinement potential while the potential $W$ is an interaction potential.  The particles are said to be in \textit{mean-field interaction} as they interact through their empirical measure $\frac{1}{N}\sum_{i=1}^N \delta_{X_t^{i,N}}$.
To this high-dimensional particle system we associate the following McKean-Vlasov nonlinear stochastic differential process in $\R^d$
\begin{equation}\label{MKV}
\begin{cases}
\mathrm{d}\bar{X}_t = \sqrt{2}\,\mathrm{d}B_t - \nabla V(\bar{X}_t)\,\mathrm{d}t - \nabla W \ast u_t(\bar{X}_t)\,\mathrm{d}t, \\
\mathcal{L}(\bar{X}_t) = u_t, \mathcal{L}(\bar{X}_0) = u_0 \\
\end{cases}
\end{equation}

where $(B_t)$ is a Brownian motion on $\R^d$ and $\mathcal{L}(X)$ denotes the law of the random variable $X$.
Such dynamics appear in the modeling of granular media, as in \cite{BCP97} and \cite{bccp98}, and more recently in machine learning theory, as in \cite{MMN18}, \cite{SS20}.

\bigskip

\textbf{Limit particle.} We denote by $(\bar{X^i})$ a solution of \eqref{MKV} when the Brownian motion $(B_t)$ is replaced by $(B_t^i)$. Moreover, we suppose that the $\bar{X}_0^i$ are independent and identically distributed, with law $u_0$. 

\bigskip
\textbf{Assumptions $\Hypp$. }We say that $V$ and $W$ satisfy assumptions $\Hypp$ if 
\begin{itemize}
    \item\label{H1} $\grad W$ is odd: $\grad W(-x) = - \grad W (x)$.
    \item\label{H2} $\left( \grad W(x) - \grad W(y)\right) \cdot \left( x-y \right) \geq \lambda_W\abs{x-y}^2 $, $\lambda_W \in \R$.
    \item\label{H3} $\left( \grad V(x) - \grad V(y) \right) \cdot \left( x-y \right) \geq \lambda_V\abs{x-y}^2$, $\lambda_V > 0$.
    \item\label{H4} $\lambda_V + 2\lambda_W > 0$.
    \item\label{H5} $\grad W$ is locally Lipschitz with polynomial growth $\abs{\grad W (x) - \grad W (y)} \leq M( \abs{x-y} \wedge 1)\left(1 + \abs{x}^m + \abs{y}^m\right)$ where $m \in \N$, $M > 0$.
    \item\label{H6}  $\grad V$ is locally Lipschitz with polynomial growth.
\end{itemize}

The constant $m$ only refers to the degree of polynomial growth of $\grad W$. While $\grad V$ is also assumed to have polynomial growth, its specific degree does not appear in our proofs and results and is therefore omitted. We do not assume convexity of the interaction potential $W$: $\lambda_W$ is not necessarily positive. This lack of convexity is compensated by the convexity of the confinement potential $V$, through the condition $\lambda_V + 2 \lambda_W > 0$.

\bigskip

\textbf{Propagation and generation of chaos.}
Given  a probability measure $\mu^N$ on a space $E^N$ we denote by $\mu^{(k,N)}$ its first $k$-point marginal.
As first introduced in \cite{Kac56}, a sequence of symmetric probability measures $(\mu^N)$ on a space $E^N$ is said to be $\nu$-chaotic, where $\nu$ is a measure on $E$, if for all $k \geq 1$, $\mu^{(k,N)}  \longrightarrow \nu^{\otimes k}$ when $N \longrightarrow + \infty$ (see \cite[Def. 3.18]{CD22b} for more details). 
We say that \textit{propagation of chaos} holds if the $u_0$-chaoticity of $u_0^{(N)}$ implies that $u_t^{(N)}$ is $u_t$-chaotic for all $t \geq 0$.
To prove \textit{propagation of chaos}, we are interested in a \textit{coupling method} as initiated by McKean and Sznitman in \cite{McKean1969}, \cite{Sznitman}. Assuming that the initial conditions of the particle system are \textit{i.i.d}, their idea is to compare the paths of the particle system to those of a system with $N$ \textit{i.i.d} processes $\bar{X}_t^i$ with common law $u_t$, with initial conditions satisfying $X_0^{i,N} = \bar{X}_0^i$. In the bounded and globally Lipschitz setting, it yields results of the form 
\begin{equation*}
    \frac{1}{N}\sum_{i=1}^N\E{\sup_{0 \leq t \leq T}\abs{X_t^{i,N} - \bar{X}_t^i}^2} \leq \frac{\Constinit}{N}e^{\Constinit t}.
\end{equation*}

 As a consequence, one obtains propagation of chaos in Wasserstein distance \eqref{defwasserstein} on finite time horizon :
\begin{equation*}
    \sup_{0 \leq t \leq T} \frac{1}{N}\mathbf{W}_2^2\left(u_t^{(N)},u_t^{\otimes N}\right) \leq \frac{\Constinit}{N}e^{\Constinit T}.
\end{equation*}
 While this result is interesting, assuming that our potentials are bounded and globally Lipschitz is restrictive, and for instance does not enable us to work with $W = \abs{x}^3/3$, which is of interest in the modeling of granular media \cite{bccp98}.
To extend this study to locally Lipschitz potentials with polynomial growth, the main tool is moment control. On a finite time horizon, asking the hessian of $V$ and $W$ to be bounded below (or a one-sided Lipschitz condition for more general drifts and diffusion coefficients) allows for moment control on $[0,T]$ and generalization of the previous result. This is the framework adopted in  \cite{BRTV98} and \cite{BCC10} in which exponential moments are used.

Another improvement would be to reach uniform-in-time bounds. In that regard, we need some dissipativity of the system, assuming that our potentials are convex (or that the one-sided Lipschitz condition holds with positive constants). This is the framework adopted by \cite{malrieu} with \textit{i.i.d} initial conditions. It is ensured by our set of assumptions $\Hypp$ through the condition $\lambda_V + 2\lambda_W > 0$.
Our first result extends \cite[Theo. 3.3]{malrieu}  to the case where the initial conditions of the particle system are neither \textit{i.i.d} nor exchangeable. Moreover, we will allow for a non-convex $W$ through the condition $\lambda_V + 2\lambda_W > 0$. 
Assuming non-\textit{i.i.d} initial conditions yields uniform-in-time  bounds for \textit{generation of chaos} in Wasserstein distance as follows:
\begin{equation}\label{typeborneW2}
    \mathbf{W}_2^2\left(u_t^{(N)},u_t^{\otimes N}\right) \leq e^{-\Const t}\mathbf{W}_2^2\left(u_0^{(N)},u_0^{\otimes N}\right) + \Constinit.
\end{equation}
The idea behind such a bound is that while the initial conditions might not be chaotic, the system will forget its initial conditions exponentially fast in time, and become chaotic. 
Note that there are other forms of generation of chaos. For instance in \cite{DMT13}, under similar assumptions but with an invariant by rotation interaction potential and under the non-chaotic regime $X_0^{i,N} = \bar{X}_0^i = X_0$, they prove generation of chaos in the sense that for all $f_1,f_2$ Lipschitz continuous functions and $\varepsilon > 0$, there exist $t_0(\varepsilon)$ and $N_0(\varepsilon)$ such that 
$
  \displaystyle  \sup_{N \geq N_0(\varepsilon)}\sup_{t \geq t_0(\varepsilon)} \mathrm{Cov}(f_1(X_t^{1,N}),f_2(X_t^{2,N})) \leq \varepsilon$.
\newline
Estimates of the form \eqref{typeborneW2} imply a convergence in $1/\sqrt{N}$ for the particle system in $\mathbf{W}_2$ distance. Similar estimates have been obtained in \cite{BDM25} with a better convergence rate of $1/N$ working, however, with a quadratic confinement potential, weak interaction and with smoother potentials $W \in \mathcal{W}^{d+3}\cap \mathrm H^s(\R^d)$ where $s > d/2+5$, which precludes potentials with polynomial growth.

\bigskip

Although the Wasserstein distance is a natural and efficient framework to quantify the propagation and generation of chaos, one can obtain finer and stronger results when turning to relative entropy estimates. This choice is motivated by the underlying dynamics of the McKean-Vlasov equation, which naturally dissipates relative entropy \eqref{defentropie}. The relative entropy controls the total variation (TV) norm through  the Csisz\'ar-Kullback-Pinsker inequality and \cite[Lemma 3.9]{DM01}, which shows that if $u_t^{(N)}$ is exchangeable
\begin{equation*}
   \NRM{u_t^{(k,N)} - u_t^{\otimes k}}_{\mathrm{TV}}^2 \leq  2\Ent{u_t^{(k,N)}}{u_t^{\otimes k}}\leq \frac{2k}{N}\Ent{u_t^{(N)}}{u_t^{\otimes N}}.
\end{equation*}

Finding bounds on $\Ent{u_t^{(N)}}{u_t^{\otimes N}}$ will be our main tool to provide propagation of chaos estimates. Note that while the previous inequality only holds for exchangeable measures, our bounds on $\Ent{u_t^{(N)}}{u_t^{\otimes N}}$ will not require exchangeability. Also note that according to \cite[Theo. 1]{OV}, convergence in relative entropy implies convergence in $\mathbf{W}_2$ distance, provided the target measure satisfies a Talagrand transportation inequality, which is ensured in our setting by the propagation of Log-Sobolev inequalities \eqref{logsob}.

Our aim will be to find bounds of the form 
\begin{equation}\label{typeborne}
    \frac{\Ent{u_t^{(N)}}{u_t^{\otimes N}}}{N} \leq \varepsilon(N)\left(\Constinit + e^{-\Constinit t}\Ent{u_0^{(N)}}{u_0^{\otimes N}}\right),
\end{equation}

where $\varepsilon(N) \longrightarrow 0 $ when $N \longrightarrow + \infty$.

This idea of \textit{propagation of chaos} was first introduced in \cite{GPV88}, and was subsequently expanded in \cite{OVa91} in dimension $1$ and in \cite{Yau91} where Yau showed that control of the entropy dissipation in time can be turned into a bound on the entropy itself for Ginzburg-Landau systems.
This idea of studying the time evolution of the entropy was later on successfully adapted in \cite{malrieu} for McKean-Vlasov equations to reach uniform-in-time entropy bounds in the general polynomial growth case, exploiting logarithmic Sobolev inequalities, and the Bakry-Emery $\Gamma_2$ criterion (\cite{BE85}, \cite{logsob}). However, while Malrieu claims that the rate $\varepsilon(N)$ can be chosen as $\varepsilon(N) = 1/N$ when the potentials exhibit polynomial growth of any order, extending the previous method to non-globally Lipschitz forces seems complicated. Indeed, the emergence of unclosed high order terms in the general polynomial growth case makes it difficult to reach such a rate using this method alone, as we will discuss in more detail later in Remark \ref{remarquemalrieu}. This entropic approach with log-Sobolev inequality in the globally Lipschitz case or bounded case has also seen recent advancements, in particular in \cite{Lac23} and 
\cite{LL23} where Lacker and Le Flem prove the sharp rate of convergence  $\Ent{u_t^{(k,N)}}{u_t^{\otimes k}} = \mathcal{O}\left(k^2/N^2\right)$ via the BBGKY hierarchy.
The entropic approach is quite robust when it comes to handling singularities. For instance, \cite{JW18} achieved fine convergence rates for some non-regular interaction kernels, such as $\grad W(x) = x/\abs{x}^{\alpha+1}, 0 < \alpha \leq 1$, and in particular for the Biot-Savart kernel, albeit on finite time horizons. Their result was adapted in \cite{GLBM24} to a uniform-in-time bound for the Biot-Savart kernel on the torus. To secure uniform-in-time bounds, time-uniform logarithmic Sobolev inequalities have been proven essential, as  studied in \cite{BG10}, \cite{malrieu} and recently generalized to a larger class of drifts in \cite{monmarche-ren-wang}. Pushing the boundaries further, \cite{RS25} establishes uniform-in-time bounds analogous to \eqref{typeborneW2} for highly singular repulsive Riesz potentials, among which $W(x) = 1/x^\alpha, \alpha < d$ and $W(x) = - \ln(\abs{x})$. The usual relative entropy dissipation method and Log-Sobolev inequalities are no longer suitable in this regime, due to the high irregularity of the interaction, but can be mimicked using a modulated free energy and logarithmic Sobolev inequalities. 

While in the previous examples the singularities are concentrated on the diagonal (e.g $ \grad W(x) = x/\abs{x}^{\alpha+1}$ with $0 <\alpha \leq 1$ in \cite{JW18}), they remain regular and Lipschitz elsewhere. We will focus on potentials that are smooth but with polynomial growth at infinity. In that regard, our main results use the classical entropic approach, coupled with time-uniform propagated log-Sobolev inequalities to secure uniform-in-time bounds. We will prove that uniform-in-time generation of chaos through bounds of the form \eqref{typeborne} can be obtained in the general polynomial growth case, with different rates $\varepsilon(N)$. These rates will depend on the growth of the interaction term, the dimension $d$ of the system and the assumptions made on the initial conditions in terms of moments, or whether they satisfy a logarithmic Sobolev inequality.

\bigskip

\textbf{Main results.} The different constants that appear throughout the statement of our results will be detailed at the end of the introduction. Our first result establishes pathwise generation of chaos through the following results:

\begin{theo}[Pathwise propagation of chaos]\label{coupling}\label{thcouplage}
Under $\Hypp$, assume that $u_0 \in \mathcal{P}_{2m^2}(\R^d)$ and \newline $u_0^{(N)} \in \mathcal{P}_2(\R^{dN})$.
Then, for all $t \geq 0$,
\begin{equation*}
    \sum_{i=1}^N \E{\abs{X_t^{i,N}- \bar{X}_t^i}^{2}} \leq 2e^{-2\beta t} \sum_{i=1}^N \E{\abs{X_0^{i,N}- \bar{X}_0^i}^{2}} + \Constinit.
\end{equation*}
Moreover, if  $u_0 \in \mathcal{P}_{2m^2p}(\R^d)$ and $u_0^{(N)} \in \mathcal{P}_{2p}(\R^{dN})$ for $p \in \N^*$ and either $d = 1$ or $W = \lambda_W \frac{\abs{x}^2}{2}$, then, for all $t \geq 0$:
\begin{equation*}
    \sum_{i=1}^N \E{\abs{X_t^{i,N}- \bar{X}_t^i}^{2p}} \leq 2^pe^{-2p\beta t} \sum_{i=1}^N \E{\abs{X_0^{i,N}- \bar{X}_0^i}^{2p}} + \frac{\Constinit}{N^{p-1}}.
\end{equation*}

\end{theo}

\begin{erem}
We do not assume that both systems start from the same initial condition; namely, we do not require that $X_0^{i,N} = \bar{X}_0^i$. Therefore, it is not necessary to have independent and identically distributed initial conditions on the particle system. In fact, no exchangeability assumption is required on the particle system.
This is reflected in our estimates since under a coupling with identical initial conditions, the term $\sum_{i=1}^N  \E{\abs{X^{i,N}_0 - \bar{X^i_0}}^{2p}}$ vanishes and the corresponding contribution disappears. It is interesting to note that over time, the system forgets its initial conditions exponentially fast.
\end{erem}
\begin{erem}
    In the second part of Theorem \ref{thcouplage}, we take advantage of higher moments to reach better convergence rates. However, we have only been able to do so in dimension $1$ and in the quadratic case.
\end{erem}

The proof we provide follows the strategy introduced in \cite{Sznitman}, \cite{BRTV98}, \cite{malrieu}, adapted to handle non-exchangeability and non-independence on the particle system, as well as non-uniform convexity.

Since $\left(X_t^{1,N},\cdots,X_t^{N,N}\right)$ and $\left(\bar{X_t^1},\cdots,\bar{X_t^N}\right)$ are random variables with respective distributions $u_t^{(N)}$ and $u_t^{\otimes N}$, we derive the following corollary:
\begin{ecor}[Propagation of chaos in Wasserstein distance]
Assume that $u_0^{(N)} \in \mathcal{P}_2(\R^{dN})$  and \newline $u_0 \in \mathcal{P}_{2m^2}(\R^d)$.
Then
\begin{equation*}
        \mathbf{W}_2^2\left(u_t^{(N)},u_t^{\otimes N}\right) \leq 2e^{-2\beta t}\mathbf{W}_2^2\left(u_0^{(N)},u_0^{\otimes N}\right) + \Constinit.
\end{equation*}

Moreover, if $u_0 \in \mathcal{P}_{2m^2p}(\R^d)$ and $u_0^{(N)} \in \mathcal{P}_{2p}(\R^{dN})$ for $p \in \N^*$ and that either  $d = 1$ or $W = \lambda_W \abs{x}^2/2$, then 

\begin{equation*}
    \mathbf{W}_{2p}^{2p}\left(u_t^{(N)},u_t^{\otimes N}\right) \leq 2^pe^{-2p \beta   t}\mathbf{W}_{2p}^{2p}(u_0^{(N)},u_0^{\otimes N}) + \frac{\Constinit}{N^{p-1}}.
\end{equation*}
\end{ecor}
\bigskip
The rest of our results are entropic generation of chaos estimates of the form \eqref{typeborne} with different rates $\varepsilon(N)$:

\begin{theo}\label{th4}
    Under $\Hypp$, assume that for all $p \in \N$, $\displaystyle \sup_{1 \leq i \leq N , N \in \N}\int \abs{x_i}^{2p}\mathrm{d}u_0^{(N)} < + \infty$. Moreover, assume that $u_0$ satisfies a logarithmic Sobolev inequality with constant $c_0$. 
    Then, for any $0 < \alpha < 1$ :
\begin{equation*}
        \frac{\Ent{u_t^{(N)}}{u_t^{\otimes N}}}{N} \leq \frac{\Constinit(\alpha)}{N^\alpha}\left(1 + e^{-\gamma t}\Ent{u_0^{(N)}}{u_0^{\otimes N}}\right).
\end{equation*}

\end{theo}

\bigskip 

In the case where the initial particles are i.i.d. and that their common law satisfies a log-Sobolev inequality, the previous result entails a  convergence rate of order $1/N^\alpha$ for all $0 < \alpha < 1$.

At the expense of an initial log-Sobolev inequality that does not depend on $N$ for the particle system, Theorem \ref{th4} can be tightened with a Markov-like argument to the following:

\begin{theo}\label{log}
  Under $\Hypp$, assume that the law of $u_0^{(N)}$ satisfies a log-Sobolev inequality with constant independent of $N$ and that $\displaystyle \sup_{1 \leq i \leq N, N \in \N} \int \abs{x_i}^{2m^2}\mathrm{d}u_0^{(N)} < +\infty$. In addition, assume that  $u_0$ satisfies a log-Sobolev inequality with constant $c_0$.
    Then, for all $t \geq 0$:
    \begin{equation*}
    \frac{\Ent{u_t^{(N)}}{u_t^{\otimes N}}}{N} \leq  \Constinit\frac{\ln(N)^m}{N}\left(1 + e^{-\gamma t}\Ent{u_0^{(N)}}{u_0^{\otimes N}}\right).
    \end{equation*}

    \end{theo}
Likewise, when the initial conditions are \textit{i.i.d}, this entails a convergence rate of $\ln(N)^m/N$.
\\
It is interesting to note that this result is the only one that requires a log-Sobolev inequality on $u_0^{(N)}$. This log-Sobolev inequality is used differently from the one on $u_0$, as it is used to propagate Gaussian moments on the particle system.
\begin{erem}\label{remarquemalrieu}
    Recall that \cite[Theo 3.3]{malrieu} claims a convergence rate of $1/N$, without providing a detailed proof. In particular, there appears to be a gap in the penultimate equation of the calculations following Lemma 3.15: bounding the stated integral is non-trivial due to the emergence of high-order moments. The sketch of the proof given in \cite{CD22b} is similarly flawed, as was discussed with the authors. Indeed, the error arises from a wrongful symmetrization in the fourth equation starting from the end of Lemma 3.4, where the term $\abs{\bar{X}_t}^{2p-1}\bar{Y}_t$ should be replaced with $\abs{\bar{Y}_t}^{2p-1}\bar{Y}_t$. Furthermore, unlike what is stated in the proof of Corollary 2.1 (page 41), it is not square moments of $\bar{X}_t^i$ that appear in the polynomial growth case, and the actual terms are the same as those that appear in Malrieu's proof as we discussed previously. The term $\ln^m(N)$ is a technical term that appears in our attempts to handle that difficulty. More details and calculations will be given in Remark \ref{remarquemoments}.
\end{erem}

We have however been able to recover Malrieu's result in the $\grad W$ Lipschitz case and when $d = 1$. It is the object of the following results:
\begin{theo}\label{cas lipschitz}
    Under $\Hypp$, assume that $\grad W$ is Lipschitz and that $u_0$ satisfies a log-Sobolev inequality with constant $c_0$. Then, for all $t \geq 0$:
    \begin{equation*}
               \frac{ \Ent{u_t^{(N)}}{u_t^{\otimes N}}}{N} \leq  \frac{\Constinit}{N}\left(e^{-\gamma t}\Ent{u_0^{(N)}}{u_0^{\otimes N}} + 1\right).
    \end{equation*}
\end{theo}

When the initial conditions of the particle system are \textit{i.i.d}, one can choose $\bar{X}_0^i = X_0^{i,N}$. In that case, $\Ent{u_0^{(N)}}{u_0^{\otimes N}} = 0$ and we recover Malrieu's result. We also recover Malrieu's result in the one dimensional case as follows:

\begin{theo}\label{thd=1}
  Under $\Hypp$, assume that  $\displaystyle \sup_{1 \leq i \leq N, N \in \N} \int \abs{x_i}^{4m^2} u_0^{(N)} < + \infty$ and that $u_0$ satisfies a log-Sobolev inequality with constant $c_0$. If moreover either $W = \lambda_W \frac{\abs{x}^2}{2}$ or $ d = 1$, then, for all $t \geq 0$:
\begin{equation*}
        \frac{\Ent{u_t^{(N)}}{u_t^{\otimes N}}}{N} \leq  \frac{\Constinit}{N} +e^{-\frac{1}{2\CLSI}t}  \frac{\Ent{u_0^{(N)}}{u_0^{\otimes N}}}{N} + e^{-\gamma t}\frac{\Constinit}{\sqrt{N}}\mathbf{W}_{4}^2\left(u_0^{(N)},u_0^{\otimes N}\right).
    \end{equation*}

  \end{theo}

  The proof of Theorem \ref{cas lipschitz} will be given in the second section as a preliminary to the proof of Theorem \ref{log}, while Theorem \ref{thd=1} will be proven at the end of the paper.
\bigskip

\textbf{Definitions and notations.} In what follows, we will denote by $\mathrm{C}$ any positive constant that depends on $V,W$ and the dimension $d$. We also denote by $\Constinit$ any positive constant that moreover depends on the initial conditions of any of the two systems. We define the constant $\beta$ by  $\beta = \lambda_V$ if $\lambda_W \geq 0$, $\lambda_V + 2 \lambda_W$ otherwise.
Given a probability measure $\mu$ we define the relative entropy of a positive $\mu$-measurable function $f$  by \begin{equation*}
    \Ent{f}{\mu} \coloneqq \int f \ln f d\mu - \int f d\mu \ln \int f d\mu 
\end{equation*}

Similarly, for $\mu$ and $\nu$ two probability measures, the relative entropy $\Ent{\mu}{\nu}$ is defined by 
        \begin{equation}\label{defentropie}
        \left\{
        \begin{aligned}
            &\int \frac{d\mu}{d\nu}\ln\left(\frac{d\mu}{d\nu}\right)d\nu \quad \text{if} \quad \mu \ll \nu \\
            &+ \infty \quad \text{otherwise.}
        \end{aligned}
        \right.
    \end{equation}
    For $p \in \N^*$ let $\mathcal{P}_p\left(\R^{dN}\right)$ be the set of probability measures on $\R^{dN}$ with a finite moment of order $p$.
The Wasserstein distance between two probability measures $\mu, \nu$ in $\mathcal{P}_{p}\left(\R^{dN}\right)$ is defined by 
\begin{equation}\label{defwasserstein}
    \mathbf{W}_p^p\left(\mu,\nu \right) \coloneqq \inf_{\pi} \int_{\R^{2dN}} \sum_{i=1}^N\abs{x^i-y^i}^p \mathrm{d}\pi(x,y),
\end{equation}
where $\abs{\cdot}$ is the Euclidean norm on $\R^d$, $x = (x^1,\cdots, x^N)$ and the infimum is taken on all couplings $\pi$ of marginals $\mu$ and $\nu$. \newline
We say that $\mu$ satisfies a log-Sobolev inequality with constant $\CLSI$ if for all suitable positive functions $f$ we have
\begin{equation}\label{logsob}
    \Ent{f}{\mu} \leq \CLSI \int \frac{\abs{\grad f}^2}{f}d\mu
\end{equation}
    The right-hand side is the Fisher information $\mathcal{I}(f|\mu)$.

If $c_0$ is a log-Sobolev constant for $u_0$, we define $\CLSI = \max(c_0,\frac{2}{\lambda_V + \lambda_W})$ and $\gamma = \min(2\beta,\frac{1}{2\CLSI})$.

\bigskip
\textbf{Organization of the paper.} The paper is organized as follows: in the first section, we prove uniform-in-time propagation of moments for both the particle system and the limit system. The proof is adapted from \cite{CGM06}, to tackle non-uniform convexity and non-exchangeability. In the second section, we establish pathwise propagation of chaos through Theorem \ref{thcouplage}. In the third section, we prove entropic propagation of chaos through Theorems \ref{th4},\ref{log}, \ref{cas lipschitz}, \ref{thd=1}.
\bigskip
\section{Moment estimates and well posedness}
In this section, we prove propagation of moments for both the particle system and the limit system. We can here work with a weaker set of assumptions than $\Hypp$. In this section we will use the following set of assumptions $\Hyp $:

\begin{itemize}
    \item\label{H_11} $\grad W$ is odd: $\grad W(-x) = - \grad W (x)$;
    \item\label{H_12} $\left( \grad W(x) - \grad W(y)\right) \cdot \left( x-y \right) \geq \lambda_W\abs{x-y}^2 - C_W  $ where $C_W \geq 0$, $\lambda_W \in \R$;
    \item\label{H_33} $\left( \grad V(x) - \grad V(y) \right) \cdot \left( x-y \right) \geq \lambda_V\abs{x-y}^2 - C_V $ where $C_V \geq 0$, $\lambda_V > 0$;
    \item\label{H_4} $\lambda_V + \lambda_W > 0$;
    \item\label{H_5} $\grad W$ is locally Lipschitz with polynomial growth:$\abs{\grad W (x) - \grad W (y)} \leq M( \abs{x-y} \wedge 1)\left(1 + \abs{x}^m + \abs{y}^m\right)$ where $m \in \N$, $M > 0$;
    \item\label{H_6}  $\grad V$ is locally Lipschitz with polynomial growth.
\end{itemize}

Allowing $C_V > 0$ and $C_W > 0$ covers non-globally convex settings, an important example being the double well potential $\abs{x}^4 - \abs{x}^2$.
\subsection{Moment estimates for the particle system}
We first establish uniform-in-time \textit{propagation of moments} for the particle system. We begin by estimating the difference between two particles. This is inspired by \cite{CGM06} and adapted to handle non-uniform convexity (\textit{i.e} a negative $\lambda_W$) and non-exchangeability on the initial laws.

\begin{eprop}\label{controleecartparticule}
    Let $1 \leq i,j \leq N, k \in \N^*$. Under $\Hyp$ assume that $\E{\abs{X_0^{i,N} - X_0^{j,N}}^{2k}} < \infty$. Then
    \begin{equation*}
        \sup_{t \geq 0} \E{\abs{X_t^{i,N} - X_t^{j,N}}^{2k}} \leq \Const(k)\left( 1 +  \E{\abs{X_0^{i,N}- X_0^{j,N}}^{2k}}\right).
    \end{equation*}
\end{eprop}

\begin{eproof}
Let $A_k(t) := \E{ \abs{X_t^i - X_t^j}^{2k}}$.  It\={o}'s formula yields 
\begin{align*}
    A_k(t) &= A_k(r) + 2k(2k + d - 2)\int_r^tA_{k-1}(s)ds \\&\quad - \frac{2k}{N}\E{\sum_{l=1}^N\int_r^t \abs{X_s^{i,N} - X_s^{j,N}}^{2k-2}\left(X_s^{i,N} - X_s^{j,N}\right)\cdot\left(\grad W(X_s^{i,N} - X_s^{l,N}) - \grad W (X_s^{j,N} - X_s^{l,N})\right) ds}\\&\quad - 2k \E{\int_r^t {\abs{X_s^{i,N} - X_s^{j,N}}^{2k-2} \left( X_s^{i,N} - X_s^{j,N}\right) \cdot \left( \grad V (X_s^{i,N}) - \grad V (X_s^{j,N})\right)\big)}ds}.
\end{align*}
Differentiating and using assumptions $\mathrm{H}_1(2)$ and $\mathrm{H}_1(3)$ yields   

\begin{equation*}
    A_k'(t) \leq 2k(2k+d-2)A_{k-1}(t) + 2k(C_V + C_W)A_{k-1}(t) - 2k(\lambda_V + \lambda_W) A_k(t).
\end{equation*}

As soon as $\lambda_V + \lambda_W > 0$, applying Grönwall's lemma and an induction yield

\begin{equation*}
    \sup_{t \geq 0}  \E{\abs{X_t^{i,N} - X_t^{j,N}}^{2k}} \leq \Const(k)\left( 1 +  \E{\abs{X_0^{i,N} - X_0^{j,N}}^{2k}}\right).
\end{equation*}
\end{eproof}

This estimate allows us to derive moment bounds.

\begin{eprop}\label{controlemomentsparticule}
        Let $1 \leq i \leq N$. Under $\Hyp$, assume that the initial laws have a finite moment of order $2mk$ and that the sequence $\left(\frac{1}{N}\sum_{j=1}^N\E{\abs{X^{i,N}_0 - X^{j,N}_0}^{2mk}}\right)_{ 1 \leq i \leq N, N \in \N}$ is bounded. 
        Then
    \begin{equation*}
       \sup_{t \geq 0, 1 \leq i \leq N,  N \in \N} \E{\abs{X_t^{i,N}}^{2k}} \leq \Constinit(k).
    \end{equation*}
\end{eprop}
\begin{erem}
    The boundedness assumption is automatically satisfied when the particle system is such that $\displaystyle \sup_{i,N}\E{\abs{X_0^{i,N}}^{2mk}} < + \infty$. If one does not suppose a uniform bound in $i$, the result still holds with a constant $\Const (k,i)$.
\end{erem}
\begin{eproof}
Let $B_k(t) := \E{\abs{X_t^{i,N}}^{2k}}$. It\={o}'s formula yields

\begin{align*}
    B_k(t) - B_k(r)&=  2k(2k + d - 2)\int_r^t B_{k-1}(s)ds - 2k \E{\int_r^t \abs{X_s^{i,N}}^{2k-2} X_s^{i,N}\cdot \grad V (X_s^{i,N})}ds \\ &\quad - \frac{2k}{N}\sum_{j=1}^N \E{\int_r^t  \abs{X_s^{i,N}}^{2k-2} X_s^{i,N} \cdot \grad W (X_s^{i,N} - X_s^{j,N})} ds \\
    &\leq 2k(2k +d -2)\int_r^t B_{k-1}(s) ds -2k\lambda_V\int_r^t B_k(s) ds + 2kC_V\int_r^tB_{k-1}(s)ds \\ &\quad -2k \int_r^t \grad V(0)\cdot \abs{X_s^{i,N}}^{2k-2}X_s^{i,N}ds -\frac{2k}{N}\sum_{j=1}^N \E{\int_r^t \abs{X_s^{i,N}}^{2k-2}X_s^{i,N} \cdot\grad W (X_s^{i,N} - X_s^{j,N})} ds \\
    &\leq 2k(2k+ d + C_V -2)\int_r^t B_{k-1}(s)ds \\ &\quad -2k(\lambda_V - \varepsilon)\int_r^tB_k(s)ds + 2kC_\varepsilon\abs{\grad V(0)}(t-r) \\ 
    & \quad -\frac{2k}{N}\sum_{j=1}^N \E{\int_r^t \abs{X_s^{i,N}}^{2k-2} X_s^{i,N} \cdot \grad W (X_s^{i,N} - X_s^{j,N})ds}
\end{align*} 
where the first inequality is due to the convexity of $V$, the second one to the fact that for all $\varepsilon > 0$ there exists $C_\varepsilon > 0$ such that $\abs{\grad V(0)}\cdot\abs{x}^{2k-1} \leq \varepsilon \abs{x}^{2k} + C_\varepsilon$.

We still need to control the interaction term. We use the polynomial growth of $\grad W$, that is $(\mathrm{H}_1(5))$ to find
\begin{equation*}
    \abs{X_s^{i,N}}^{2k-2}\abs{ \grad W (X_s^{i,N} - X_s^{j,N})\cdot X_s^{i,N}}\leq \abs{X_s^{i,N}}^{2k-1} M\left( 1 + \abs{X_s^{i,N} - X_s^{j,N}}^m \right).
\end{equation*}

Hölder's inequality yields  

\begin{equation*}
    \E{\abs{X_s^{i,N}}^{2k-2}\abs{ \grad W (X_s^{i,N} - X_s^{j,N})\cdot X_s^{i,N}} }\leq M\left(\E{\abs{X_s^{i,N}}^{2k}}\right)^{\frac{2k-1}{2k}}\left( \E{\left( 1 + \abs{X_s^{i,N} - X_s^{j,N}}^m \right)^{2k}}\right)^{\frac{1}{2k}}.
\end{equation*}

Now, averaging over $j$ only affects the second term of the product, which is bounded in $i,N$ and in time, thanks to Jensen's inequality, Proposition \ref{controleecartparticule} and our assumption. Therefore, the averaged term can be bounded by $\Const B_k^{\frac{2k-1}{2k}}$. Hence  
\begin{equation*}
    B_k(t) - B_k(r) \leq 2k(2k+d+ C_V-2)\int_r^t B_{k-1}(s)ds +2kC_\varepsilon(t-r) - 2k(\lambda_V - \varepsilon) \int_r^t B_k(s) ds + 2kC\int_r^t B_{k}^{\frac{2k-1}{2k}}(s)ds.
\end{equation*}
Differentiating gives
\begin{equation*}
    B_k'(t) \leq 2k(2k+d + C_V-2)B_{k-1}(t) + 2kC_\varepsilon - 2k(\lambda_V-\varepsilon) B_k(t) + 2k\Const\left(B_{k}(t)\right)^{\frac{2k-1}{2k}}.
\end{equation*}

To conclude, we use the fact that for all $\varepsilon > 0$, there exists $c > 0$ such that  for all $a > 0$, $Ca^{\frac{2k-1}{2k}} \leq c_\varepsilon + \varepsilon a$. This yields $B'_{k}(t) \leq 2k(C_\varepsilon + c_\varepsilon) - 2k(\lambda_V - 2\varepsilon)B_k + 2k(2k+d-2 + C_V)B_{k-1}$. We choose $\varepsilon < \frac{\lambda_V}{2}$ and the conclusion follows by Grönwall's lemma and an induction.
\end{eproof}

\subsection{Moment estimates for the limit particle}

We now apply the same strategy to the limit particle. This is again adapted from \cite{CGM06}, which deals with the square moments only. We extend their result to every moments of the limit particle. 
To control the interaction term, we begin by controlling the gap between two independent limit particles. 

\begin{eprop}\label{ecartdeuxpartlim}
        Let $(\bar{X}_t)_{t \geq 0}$ and $(\bar{Y}_t)_{t \geq 0}$ be independent solutions of $\eqref{MKV}$. Under $\Hyp$, assume that the initial laws have a finite moment of order $2k$. Then
            \begin{equation*}
        \sup_{t \geq 0} \E{\abs{\bar{X}_t - \bar{Y}_t}^{2k}} \leq \Const(k)\left(1 + \E{\abs{\bar{X_0} - \bar{Y_0}}^{2k}}\right).
    \end{equation*}
\end{eprop}

    \begin{eproof}
        
 Let $y(t) := \E{\abs{\bar{X}_t-\bar{Y}_t}^{2k}}$.
    Then, applying It\={o}'s formula yields 

    \begin{align*}
        y'(t) &= -2k\E{\abs{\bar{X}_t - \bar{Y}_t}^{2k-2}\left(\bar{X}_t - \bar{Y}_t\right)\cdot\left( \grad V (\bar{X}_t) - \grad V (\bar{Y}_t)\right)}  \\ &\quad -2k\E{\abs{\bar{X}_t - \bar{Y}_t}^{2k-2}\left(\bar{X}_t - \bar{Y}_t\right)\cdot\left( \grad W \ast u_t(\bar{X}_t) - \grad W \ast u_t(\bar{Y}_t)\right)} \\
        &\quad+ 2k(2k+d-2)\E{\abs{\bar{X}_t - \bar{Y}_t}^{2k-2}}
    \end{align*}

    Using convexity on $V$ and $W$ gives us
    \begin{align}
        y'(t) \leq -2k(\lambda_V + \lambda_W)y(t) + 2k(C_V + C_W + 2k +d - 2)\E{\abs{\bar{X}_t - \bar{Y}_t}^{2k-2}}.
    \end{align}

    The conclusion follows by Grönwall's lemma and an induction.
    \end{eproof}
We obtain a result analogous to that of the particle system:
\begin{eprop}\label{particulelimiteborne}
       Let $k \in \N$. Let $\bar{X}$ be a solution of $\eqref{MKV}$ under $\Hyp$. Assume that $\E{\abs{\bar{X}_0}^{2mk}} < + \infty$. Then \begin{equation*}
        \sup_{t \geq 0}\E{\abs{\bar{X}_t}^{2k}} < \Const(k)\left(1 + \E{\abs{\bar{X_0}}^{2mk}}\right).
    \end{equation*} 

\end{eprop}

\begin{eproof}
        Let $y(t) := \E{\abs{\bar{X}_t}^{2k}}$.
    Then 
    \begin{equation}\label{preumdeux}
        \begin{aligned}
        y'(t) &= -2k\E{\abs{\bar{X}_t}^{2k-2}\bar{X}_t\cdot\grad V (\bar{X}_t)}-2k\E{\abs{\bar{X}_t}^{2k-2}\bar{X}_t\cdot\grad W \ast u_t (\bar{X}_t)}  + 2k(2k+d-2)\E{\abs{\bar{X}_t}^{2k-2}}\\
        &= -2k\E{\abs{\bar{X}_t}^{2k-2}\bar{X}_t\cdot\left(\grad V (\bar{X}_t) - \grad V(0) + \grad V(0) \right)} -2k\E{\abs{\bar{X}_t}^{2k-2}\bar{X}_t\cdot\grad W \ast u_t (\bar{X}_t)}
    \end{aligned}
    \end{equation}

As in the proof of Proposition \ref{controlemomentsparticule}, using  the inequality $\abs{\grad V(0)}\abs{x}^{1-\frac{1}{2k}} \leq \varepsilon \abs{x}^{2k}+ C_\varepsilon$, the first term in \eqref{preumdeux} is bounded by
\begin{equation*}
- 2k (\lambda_V - \varepsilon)y(t) + 2k(C_\varepsilon + C_V) + 2k(2k+d-2)\E{\abs{\bar{X}_t}^{2k}}. 
    \end{equation*}

Moreover,    
    \begin{equation*}
        \E{\abs{\grad W \ast u_t(\bar{X}_t)}^{2k}} \leq \E{\int \abs{\grad W (X_t - y)}^{2k}}du_t(y) 
        \leq M \E{\int \abs{1 + \abs{\bar{X}_t - y}^{2mk}}du_t(y)}
    \end{equation*}

    is bounded in time by Proposition \ref{ecartdeuxpartlim}. 

        By Hölder's inequality, the second term in \eqref{preumdeux} is bounded by: 
\begin{equation*}
     2k y(t)^{1-\frac{1}{2k}}\E{\abs{\grad W \ast u_t(\bar{X}_t)^{2k}}}^{\frac{1}{2k}} \leq C2ky(t)^{1-\frac{1}{2k}}.
\end{equation*}
We can now conclude as in Proposition \ref{controlemomentsparticule}.
\end{eproof}

\subsection{Well-posedness} 

Strong well-posedness of the particle system equation \eqref{systemedeparticules} is ensured as soon as the initial law $u_0^{(N)}$ has a second moment thanks to the convexity of $V$ and $W$ as proven in \cite[Theo. 10.2.2, p. 255]{SV06}.
It is, however, not obvious that the nonlinear limit equation \eqref{MKV} is well-posed. In the Lipschitz case, existence and uniqueness follow from a standard fixed-point argument in $\left(C\left(\left[0,T\right],\mathcal{P}_2\left(\R^d\right)\right), \sup\limits_{0 \leq t \leq T}\mathbf{W}_2\left(\cdot,\cdot\right)\right)$ as in \cite{Sznitman} or \cite{PECDR}. 
However, the same procedure  does not readily extend to the setting of polynomial growth, as moments of order $2m$ appear in the fixed-point procedure and must be bounded using similar methods as those we used in this section.

A first approach, followed in \cite{CGM06} and \cite{Meleard-1996}, consists in proving that the particle system converges through moment estimates and a tightness criterion and showing that the limit is a solution of the nonlinear martingale problem. This procedure requires moments of order $2m^2$ on the initial conditions and that the initial conditions be chaotic and exchangeable. If one is interested in generation of chaos rather than propagation of chaos, assuming chaos on the initial condition is precluded.

A fixed-point argument can  still be applied, however not in a measure space as in the Lipschitz case. The appropriate space to use a fixed-point theorem here is the function space \begin{equation*}
    \displaystyle \Lambda_T = \left\{ b : [0,T] \longrightarrow \R^d, x \mapsto b(x,t) \text{ is loc. Lip.}, \left(b(t,x) - b(t,y) \cdot x-y\right) \geq \lambda_W \abs{x-y}^2 -C_W, ||b||_T < \infty \right\},
\end{equation*}
where $||b||_T = \sup_{t \in [0,T]}\sup_{x \in \R^d} \frac{b(t,x)}{1 + ||x||^{2q}} $ and $q$ is such that $2q \geq m +1$.
 This is done in dimension 1 in \cite{BRTV98}, and adapted in the multidimensional case in \cite{Herrmann-Imkeller-Peithmann-2008}, under however slightly different assumptions. Combining the two approaches yields the following:

\begin{eprop}
    Under $\Hypp$, equation \eqref{MKV} has a unique strong solution when one of the following holds: 
    \begin{itemize}
        \item $u_0 \in \mathcal{P}_{2}(\R^d)$ if $\grad W$ is Lipschitz
        \item  $ u_0 \in \mathcal{P}_{2m^2}(\R^d)$ otherwise.
    \end{itemize}
\end{eprop}

\section{Pathwise Propagation of Chaos: Proof of Theorem \ref{thcouplage}}
In this section, we are interested in the pathwise convergence of the particle system towards the limit particle. 
For this section, all of our results will be under $\Hypp$.
We denote by $\bar{X}^i$ a solution of \eqref{MKV} driven by $B^i$. This is a coupling procedure \textit{à la Sznitman}. This section is adapted from \cite{BRTV98}, \cite{Sznitman},\cite{malrieu} which studies the case $p = 1$ to handle non-exchangeability, non-uniformly convex potential and non-independence on the initial conditions of the particle system. We also give sharp assumptions on the initial conditions, as we do not need a log-Sobolev inequality in that section, nor moments of all orders. 
\bigskip
\begin{eproof}[Theorem \ref{thcouplage}]
    Let $\Delta^i := X^{i,N}-\bar{X^i}$ and $S_t := \sum_{i=1}^N \E{\abs{\Delta_t^i}^{2p}}$.

 It\={o}'s formula yields
\begin{equation*}
\begin{aligned}\label{termefin}
    \sum_{i=1}^N \abs{\Delta_t^i}^{2p} - \sum_{i=1}^N \abs{\Delta_s^i}^{2p} &= -2p\sum_{i=1}^N\int_s^t \abs{\Delta_r^i}^{2p-2}\left(\Delta^i \cdot \grad V(X_r^i) - \grad V (\bar{X_r^{i,N}})\right)dr \\ &-\frac{2p}{N}\sum_{1 \leq i,j \leq N}\int_s^t \abs{\Delta_r^i}^{2p-2}\Delta_r^i \cdot \left(\grad W (X_r^{i,N} - X_r^{j,N}) - \grad W \ast u_r (\bar{X}_r^i)\right) dr.
\end{aligned}
\end{equation*}

\textbf{Step 1.}
By convexity of $V$ the integrand of the first term in the right-hand side is bounded by $-2p\lambda_VS_r$.

\textbf{Step 2.}
We write the integrand of the second term as 
\begin{equation}\label{tempsr}
\begin{aligned}
&-\frac{2p}{N}\sum_{1 \leq i,j \leq N}  \abs{\Delta^i_r}^{2p-2}\Delta^i_r \cdot \left(\grad W (X_r^{i,N} - X_r^{j,N}) - \grad W (\bar{X}_r^i - \bar{X}_r^j) \right)  \\ &-\frac{2p}{N}\sum_{1 \leq i,j \leq N}  \abs{\Delta_r^i}^{2p-2}\Delta_r^i \cdot \left( \grad W (\bar{X}_r^i - \bar{X}_r^j) - \grad W \ast u_r(\bar{X}_r^i) \right).
\end{aligned}
\end{equation}

\textbf{Step 2.1.} We begin with the second term in \eqref{tempsr}. Let $1 \leq i \leq N$ and consider the sum over $j$.
Hölder's inequality shows that it is bounded by  
\begin{equation}\label{etape2.1}
 \frac{2p}{N}\left(\E{\abs{\Delta_r^i}^{2p}}\right)^{\frac{2p-1}{2p}}\quad\E{\abs{\sum_{j=1}^N \grad W (\bar{X}_r^i - \bar{X}_r^j) - \grad W \ast u_r(\bar{X}_r^i)}^{2p}}^\frac{1}{2p}.
\end{equation}

We begin with the sum on $j \neq i$. Developing the exponent $2p$ yields
\begin{equation*}
    \E{\abs{\sum_{\substack{j=1 \\j \neq i}}^N \grad W (\bar{X}_r^i - \bar{X}_r^j) - \grad W \ast u_r(\bar{X}_r^i)}^{2p}} = \sum_{\ell=1}^{2p} \sum_{\substack{1 \leq j_1 \neq j_2 \neq...j_\ell \leq N \\ j_l \neq i}} \E{\prod_{m_1 + \cdots + m_\ell = 2p} a_{j_k}^{m_k}}
\end{equation*}
where $a_j = \grad W (\bar{X}_r^i - \bar{X}_r^j) - \grad W \ast u_r(\bar{X}_r^i)$.

If an exponent satisfies $m_k = 1$, then by independence conditionally on $ \bar{X}_r^i$ and the fact that $\E{a_j | \bar{X}_r^i} = 0$ for $j \neq i$ then
\begin{equation*}
    \E{\prod_{m_1 + \cdots + m_\ell = 2p} a_{j_k}^{m_k}} = 0.
\end{equation*}

This yields
\begin{equation}\label{avantccontroleprod}
    \E{\abs{\sum_{\substack{j=1 \\ j \neq i }}^N \grad W (\bar{X}_r^i - \bar{X}_r^j) - \grad W \ast u_t(\bar{X}_r^i)}^{2p}} = \sum_{\ell=1}^{p} \sum_{\substack{1 \leq j_1 \neq j_2 \neq...j_l \leq N \\ j_l \neq i} } \E{\prod_{\substack{m_1 + \cdots + m_\ell = 2p \\ m_k \geq 2 }} a_{j_k}^{m_k}}
\end{equation}

where we changed the index on the first sum because all terms with $r > p$ force an exponent $m_k $ to be equal to $1$.
Now, to find an upper bound on the expectation in \eqref{avantccontroleprod}, with $q_k = \frac{2p}{m_k}$, then $\sum_{k = 1}^r \frac{1}{q_k} = 1$ and Hölder's inequality yields

\begin{equation*}
    \E{\prod_{\substack{m_1 + \cdots + m_\ell = 2p \\ m_l \geq 2 }} a_{j_k}^{m_k}} \leq \prod_{k = 1}^r \E{\abs{a_{j_k}}^{m_kq_k}}^{\frac{1}{q_k}} = \prod_{k = 1}^r \E{\abs{a_{j_k}}^{2p}}^{\frac{m_k}{2p}} \leq \Constinit
\end{equation*}

where the last inequality is given by our moment control through Proposition \ref{particulelimiteborne} applied with  $ k = p$ after using our polynomial growth hypothesis. This entails

    \begin{equation*}
    \E{\abs{\sum_{\substack{j=1 \\ j \neq i}}^N \grad W (\bar{X_r^i} - \bar{X_r^j}) - \grad W \ast u_t(\bar{X_r^i})}^{2p}} \leq \sum_{r=1}^p\sum_{\substack{1\leq j_1 \neq j_2 \neq...j_r \leq N \\ j_\ell \neq i}  } \Constinit 
     = \Constinit \sum_{r=1}^p \binom{N}{r} \leq \Constinit N^p.
\end{equation*}

Now, we must add the term $j=i$. Using that $\E{\abs{X+Y}^{2p}}^{\frac{1}{2p}} \leq \E{\abs{X}^{2p}}^{\frac{1}{2p}} + \E{\abs{Y}^{2p}}^{\frac{1}{2p}}$ with $\abs{X}^{2p}$ being the previous sum and $Y = \grad W \ast u_r(\bar{X}_r^i)$ and using polynomial growth $(\mathrm{H}_2(5))$ and moment control on $\E{\abs{Y}^{2p}}$, which is ensured through Proposition \ref{particulelimiteborne} with $k = 2mp$ yields
\begin{equation}\label{etape2.1bis}
    \E{\abs{\sum_{j=1}^N \grad W (\bar{X}_r^i - \bar{X}_r^j) - \grad W \ast u_r(\bar{X}_r^i)}^{2p}}^{\frac{1}{2p}} = \E{\abs{X+Y}^{2p}}^{\frac{1}{2p}} \leq \Constinit\sqrt{N} + \Constinit \leq \Constinit\sqrt{N}.
\end{equation}

 Finally, injecting \eqref{etape2.1bis} in \eqref{etape2.1} yields

 \begin{equation*}
 -\frac{2p}{N}\E{\sum_{j=1}^N  \abs{\Delta_r^i}^{2p-2}\Delta_r^i \cdot \left(\grad W (\bar{X}_r^i - \bar{X}_r^j) - \grad W \ast u_r(\bar{X}_r^i) \right)} \leq C\frac{\left(\E{\abs{\Delta_r^i}^{2p}}\right)^{\frac{2p-1}{2p}}}{\sqrt{N}}
\end{equation*}

\bigskip

\textbf{Step 2.2}
Now, we study the integrand of the first term in \eqref{tempsr}
\begin{equation}\label{termecomplexe}
    -\frac{2p}{N}\sum_{i,j}  \abs{\Delta_r^i}^{2p-2}\Delta_r^i \cdot \left(\grad W (X_r^{i,N} - X_r^{j,N}) - \grad W (\bar{X}_r^i - \bar{X}_r^j) \right) = -\frac{2p}{N} \sum_{1 \leq i,j \leq N}p_{i,j}.
\end{equation}
Since $\grad W(-z) = - \grad W (z)$
\begin{equation*}
 p_{i,j} + p_{j,i} = \left(\abs{\Delta_r^i}^{2p-2}\Delta^i - \abs{\Delta_r^j}^{2p-2}\Delta_r^j \right) \cdot \left(\grad W(X_r^{i,N} - X_r^{j,N}) - \grad W(\bar{X}_r^i- \bar{X}_r^j) \right).
\end{equation*}
Note that $p_{i,i} = 0$.

\textbf{Step 2.2.1 } We begin with the case $p = 1$. Then, after symmetrization

\begin{equation*}
    p_{i,j} + p_{j,i} = (\Delta_r^i - \Delta_r^j)\cdot (\grad W (X_r^{i,N} - X_r^{j,N}) - \grad W (\bar{X}_r^i - \bar{X}_r^j)) \geq \lambda_W \abs{\Delta_r^i - \Delta_r^j}^2
\end{equation*}

If $\lambda_W \geq 0$, this is a positive term we can neglect. If $\lambda_W \leq 0$ we have the following

\begin{equation*}
    p_{i,j} + p_{j,i} \geq 2 \lambda_W(\abs{\Delta_r^i}^2 + \abs{\Delta_r^j}^2),
\end{equation*}

which  yields

\begin{equation*}
    -\frac{2}{N} \E{\sum_{1 \leq i,j \leq N} p_{i,j}} \leq -4\lambda_W \sum_{i=1}^N \E{\abs{\Delta_r^i}^{2}} = -4\lambda_W S_r.
\end{equation*}

\textbf{Step 2.2.2 } If  $ d =1$, then $x^{2p-1}-y^{2p-1}$ has the same sign as $x-y$. Therefore, by convexity of $W$

\begin{equation}\label{symmetrisationpositif}
\begin{aligned}
    p_{i,j} + p_{j,i} &= \frac{\left(\abs{\Delta_r^i}^{2p-2}\Delta_r^i - \abs{\Delta_r^j}^{2p-2}\Delta_r^j \right)}{\Delta_r^i - \Delta_r^j}\left( \Delta_r^i - \Delta_r^j \right) \cdot \left(\grad W(X_r^{i,N} - X_r^{j,N}) - \grad W(\bar{X_r^i}- \bar{X_r^j}) \right) \\
    &\geq \lambda_W \frac{\left(\abs{\Delta^i}^{2p-2}\Delta^i - \abs{\Delta^j}^{2p-2}\Delta^j \right)}{\Delta^i - \Delta^j} \abs{\Delta^i - \Delta^j}^2.
    \end{aligned}
\end{equation}

If $\lambda_W \geq 0$, this is a positive term which we can neglect. 

In the negative $\lambda_W$ case, we write \eqref{symmetrisationpositif} as
\begin{equation*}
\begin{aligned}
        p_{i,j} + p_{j,i} &\geq \lambda_W \left(\abs{\Delta_r^i}^{2p-2}\Delta_r^i - \abs{\Delta_r^j}^{2p-2}\Delta_r^j \right)\left(\Delta_r^i - \Delta_r^j\right)\\
    &= \lambda_W\left( \abs{\Delta_r^i}^{2p} + \abs{\Delta_r^j}^{2p} - \abs{\Delta_r^i}^{2p-2}\Delta_r^i\Delta_r^j - \abs{\Delta_r^j}^{2p-2}\Delta_r^i\Delta_r^j\right).
\end{aligned}
\end{equation*}Hölder's inequality on the crossed term yields, for $i \neq j$:
    
\begin{equation*}
   \E{ p_{i,j} + p_{j,i}} \geq  \lambda_W\left( \E{\abs{\Delta_r^i}^{2p}} + \E{\abs{\Delta_r^j}^{2p}} + \E{\abs{\Delta_r^i}^{2p}}^{\frac{2p-1}{2p}}\E{\abs{\Delta_r^j}^{2p}}^{\frac{1}{2p}} + \E{\abs{\Delta_r^j}^{2p}}^{\frac{2p-1}{2p}}\E{\abs{\Delta_r^i}^{2p}}^{\frac{1}{2p}}\right).
\end{equation*}

Summing on $1 \leq i,j \leq N$ yields
\begin{equation*}
    -\frac{2p}{N} \E{\sum_{1 \leq i,j \leq N}p_{i,j}} \leq - 2p\lambda_W \sum_{i=1}^N \E{\abs{\Delta_r^i}^{2p}}  -2p \lambda_W \sum_{1 \leq i,j \leq N}\left(\E{\abs{\Delta_r^i}^{2p}}\right)^{\frac{2p-1}{2p}}\left(\E{\abs{\Delta_r^j}^{2p}}\right)^{\frac{1}{2p}}
\end{equation*}

Now, we use the fact that for all positive $(a_i)_{1 \leq i \leq N}$, $\sum_{i=1}^N \left(a_i^{\frac{2p-1}{2p}}\right) \leq N^{\frac{1}{2p}}\left(\sum_{i=1}^N a_i\right)^{\frac{2p-1}{2p}}$ and that  \newline $\sum_{i=1}^N a_i^{\frac{1}{2p}} \leq N^{1-\frac{1}{2p}}\left(\sum_{i=1}^N a_i\right)^{\frac{1}{2p}}$ to prove

\begin{equation*}
    -\frac{2p}{N}\E{\sum_{1 \leq i,j \leq N} p_{i,j}} \leq -4p\lambda_W \sum_{i=1}^N \E{\abs{\Delta_r^i}^{2p}} = -4p \lambda_W S_r.
\end{equation*}

If $d \geq 2$, the same procedure does not readily work, as  $\abs{\Delta_r^i}^{2p-2}\Delta_r^i - \abs{\Delta_r^j}^{2p-2}\Delta_r^j $ is not necessarily collinear to $\Delta_r^i - \Delta_r^j$ and we cannot use convexity on this term.

\textbf{Step 2.2.3 }The case $W = \lambda_W \frac{\abs{x}^2}{2}$:

Here, the term becomes 

\begin{equation*}
    \lambda_W \left(\abs{(\Delta_r^i}^{2p-2}\Delta_r^i - \abs{\Delta_r^j}^{2p-2}\Delta_r^j )\cdot (\Delta_r^i - \Delta_r^j)\right).
\end{equation*}

As soon as $\lambda_W \geq 0$, it is a positive term, as $x \mapsto |x|^{2p-2}x$ is monotonic. 

If $\lambda_W < 0$, to derive a lower bound on this term we must get an upper bound on  

\begin{equation*}
    \E{\left(\abs{\Delta_r^i}^{2p-2}\Delta_r^i - \abs{\Delta_r^j}^{2p-2}\Delta_r^j )\cdot (\Delta_r^i - \Delta_r^j)\right)} = \E{\abs{\Delta_r^i}^{2p} + \abs{\Delta_r^j}^{2p} - \abs{\Delta_r^i}^{2p-2}\Delta_r^i \cdot \Delta_r^j - \abs{\Delta_r^j}^{2p-2}\Delta_r^j \cdot \Delta_r^i}
\end{equation*}

The same calculations as in step 2.2.2 yield the same result.

\textbf{Step 3 }
Plugging Steps 2.1 and 2.2 in \eqref{tempsr} and using Step 1 shows that 
\begin{equation*}
    S_t-S_s \leq -2p\beta\int_s^t S_r dr + \frac{\Constinit}{\sqrt{N}}\sum_{i=1}^N\int_s^t\left(\E{\abs{\Delta_r^i}^{2p}}\right)^{\frac{2p-1}{2p}}dr.
\end{equation*}

Differentiating and using Jensen's inequality with the concave function $x \mapsto x^{\frac{2p-1}{2p}}$ yields

\begin{equation*}
    S_t' \leq -2p\beta S_t + \frac{\Constinit}{\sqrt{N}}\sum_{i=1}^N \left(\E{\abs{\Delta_t^i}^{2p}}\right)^{\frac{2p-1}{2p}} \leq -2p \beta S_t + \frac{\Constinit}{N^{\frac{1}{2}-\frac{1}{2p}}}S_t^{\frac{2p-1}{2p}}.
\end{equation*}

Grönwall's lemma shows 

\begin{equation}\label{wass}
    S_t \leq \left(e^{-\beta t}S(0)^{\frac{1}{2p}} + \frac{\Constinit}{N^{\frac{1}{2}- \frac{1}{2p}}}\right)^{2p} 
\end{equation}

which finally yields

\begin{equation*}
    \sum_{i=1}^N \E{\abs{X_t^{i,N}- \bar{X}_t^i}^{2p}} \leq 2^p e^{-2p\beta t} \sum_{i=1}^N \E{\abs{X_0^{i,N}- \bar{X}_0^i}^{2p}} + \frac{\Constinit}{N^{p-1}}.
\end{equation*}
\end{eproof}

\section{Entropic Propagation of Chaos: Proofs of Theorems \ref{th4}, \ref{log}, \ref{cas lipschitz}, \ref{thd=1} }

\subsection{Entropy dynamics and propagation of log-Sobolev inequalities}
The proof of Theorems \ref{th4}, \ref{log}, \ref{cas lipschitz}, \ref{thd=1} is based on the evolution of  $H(t) := \Ent{u_t^{(N)}}{u_t^{\otimes N}}$, where $u_t^{(N)}$ is the law at time $t$ of the particle system and $(u_t)$ the law at time $t$ of the limit system.

 We have the following Fokker-Planck equation for the law of the solutions $(u_t)$ of \eqref{MKV} and $(u_t^{(N)}) $ of \eqref{systemedeparticules}:
\begin{eprop}[Fokker-Planck dynamics]\label{dynamiques}
    Let $u_0^{(N)}$ and $u_0$ be some probability measures, respectively on $\R^{dN}$ and $\R^d$. Under $\Hypp$, assume that $u_0^{(N)}$ has a finite second moment and $u_0$ a finite $2m^2$ moment.
Then, $u_t^{(N)}$ is the solution of
\begin{equation*}
         \partial_t u_t^{(N)} = \lap u_t^{(N)} + \grad \cdot \left(u_t^{(N)}\grad \psi^{(N)}\right) 
        \end{equation*}
where $\psi^{(N)}$ is defined for $x_1,...,x_N$ in $\R^d$ by
\begin{equation*}
    \psi^{(N)}(x_1,...,x_N) := \sum_{i=1}^N V(x_i) + \frac{1}{2N}\sum_{i,j = 1}^N W(x_i-x_j).
\end{equation*}
Moreover, $u_t$ is the solution of the following equation:
\begin{equation}\label{eqlim}
       \partial_t u_t = \lap u_t + \grad \cdot \left(u_t(\grad V + \grad W \ast u_t)\right)
\end{equation}

\end{eprop}

The main tool to study the evolution of $H(t) = \Ent{u_t^{(N)}}{u_t^{\otimes N}}$ will be the propagation in time of log-Sobolev inequalities for both systems.

\begin{eprop}[Time-propagation of a log-Sobolev inequality for the particle system]\label{proplsiparticle}
    Assume that $u_0^{(N)}$  satisfies a log-Sobolev inequality with constant $c_0$. Then, $u_t^{(N)}$ satisfies a log-Sobolev inequality with constant $c_0e^{-2\rho t} + \frac{2}{\rho}\left(1 - e^{-2 \rho t}\right)$ with $\rho = \lambda_V $ if $\lambda_W > 0$ and $\lambda_V + \lambda_W $ otherwise.
\end{eprop}

    \begin{eproof}
    We have to show that $u_t^{(N)}$ satisfies the $CD(\rho,\infty)$ criterion (see \cite{logsob} and \cite{malrieu} for more details)

Its infinitesimal generator is $\mathbf{L}^N F = \Delta F -\nabla\psi^{(N)} \nabla F$ where $\psi^{(N)}$ was defined in Proposition \ref{dynamiques}.

To apply the Bakry-Emery $\Gamma_2$ criterion, we must show that $\grad^2 \psi^{(N)} \geq \rho I_{dN} $.

Let $h = (h_1,\cdots,h_N) \in \R^{dN}$, $x = (x_1,\cdots,x_N) \in \R^{dN}$. Then

\begin{equation*}
    \begin{aligned}
        \grad^2 \psi^{(N)}(x) (h,h) &= \sum_{i=1}^N  \grad^2 V(x_i)h_i \cdot h_i +\frac{1}{2N}\sum_{i,j=1}^N \grad^2 W (x_i - x_j)(h_i-h_j) \cdot (h_i-h_j) \\ &\geq \lambda_V \abs{h}^2 + \frac{\lambda_W}{2N}\sum_{i,j = 1}^N \abs{h_i - h_j}^2 \\ 
        &= \lambda_V \abs{h}^2 +\lambda_W \abs{h}^2 - \frac{\lambda_W}{N}\abs{\sum_{i=1}^N h_i}^2. \\
    \end{aligned}
\end{equation*}
If $\lambda_W \geq 0 $, with Cauchy-Schwarz inequality
\begin{equation*}
        \grad^2 \psi^{(N)}(x)(h,h) \geq \lambda_V \abs{h}^2 + \lambda_W\abs{h}^2 - \lambda_W \sum_{i=1}^N\abs{h_i}^2 = \lambda_V\abs{h}^2.
\end{equation*}
If $\lambda_W < 0$, then we can neglect the last term 
\begin{equation*}
    \grad^2 \psi^{(N)}(x)(h,h) \geq (\lambda_V + \lambda_W)\abs{h}^2,
\end{equation*}
which concludes.
\end{eproof}

We now prove that the log-Sobolev inequality is propagated in time by \eqref{eqlim}.

\begin{eprop}[Time-propagation of a log-Sobolev inequality for the limit particle]\label{propaglsilimit}
Let $(u_t)$ be a solution of \eqref{eqlim}. We suppose that $u_0$ satisfies a log-Sobolev inequality with constant $c_0$. Then, $u_t$ satisfies a log-Sobolev inequality with constant 
$$ c_t = c_0e^{-2\rho t} +  \frac{2}{\rho}\left(1-e^{-2\rho t}\right),$$
where $\rho = \lambda_V + \lambda_W$.

\end{eprop}
\begin{eproof}
There are two ways to prove this result, which are both detailed in \cite{malrieu}. The first one is to use a particle system, with initial conditions $X_0^{i,N}$ independent and identically distributed with law $u_0$. Using Proposition \ref{proplsiparticle} and the convergence of $u_t^{1,N}$ to $u_t$ in $\mathbf{W}_2$ through Theorem \ref{thcouplage} concludes, as our log-Sobolev inequality is uniform in $N$. However, in the $\lambda_W \geq 0$ case, we only have the result for $\rho = \lambda_V$.

To conclude without the help of a particle system, we can prove that the infinitesimal generator satisfies the $CD(\rho,\infty)$ criterion, which is inhomogeneous this time. For more detail, one can read \cite[Theorem 23]{BG10}, \cite[Theorem 3.10]{colletmalrieu08}, \cite{monmarche-ren-wang}.
     In that case, we must show that for all suitable functions $f$, $\Gamma_2^s(f)\geq \rho \Gamma(f)$, where: 
\begin{equation*}
           \Gamma(f) := \abs{\grad f}^2, \Gamma_2^s(f) := \NRM{\grad^2 f}^2 + \left(\grad f\right)^T \grad^2 V \left(\grad f\right) + \left(\grad f\right)^T(\grad^2 W \ast u_s) \left(\grad f\right),
\end{equation*}

Then, \begin{equation*}
    \Gamma_2^s(f) \geq \left(\grad f\right)^T \grad^2 V \left(\grad f\right) + \left(\grad f\right)^T(\grad^2 W \ast u_s) \left(\grad f\right) \geq \left(\lambda_V + \lambda_W\right)\abs{\grad f}^2 = \rho\Gamma(f).
\end{equation*}
\end{eproof}
\begin{erem}
The influence of the initial constant $c_0$ decays exponentially in time.\end{erem}
\begin{erem}\label{remark30}
The constant $c_t$ \textit{a priori} depends on $t$. However, it is a convex combination of $c_0$ and $\frac{2}{\rho}$, it can be bounded by a constant that does not depend on time, for instance by $\CLSI$ which was defined as $\CLSI = \max\left(c_0,\frac{2}{\rho}\right)$.
\end{erem}
\begin{erem}
    If the initial law is a Dirac measure, then $c_0 = 0$.
\end{erem}

The following result about the relative entropy of two solutions with different vector fields is often used, we have however not found it clearly stated and proven in the literature. For completeness, we provide a proof and detailed calculations :
\begin{eprop}\label{dynamiquentropie}
    
    Let $h = h(x,t), g = g(x,t), x \in \R^d, t \geq 0$ be densities such that $\partial_t h  = \lap h + \diver(hB)$, $\partial_t g = \lap g + \diver (gA)$ with vector fields $A,B : \R^d \times \R_+ \longrightarrow \R^d$. Then
    $$\partial_t \Ent{h_t}{g_t} = -\mathcal{I}(h | g) - \int h\grad(\ln\frac{h}{g})(B-A)$$
    where $\mathcal{I}(h | g) = \int h \abs{\grad \ln{\frac{h}{g}}}^2$ is the Fisher information.
\end{eprop}
\begin{eproof}

\begin{equation*}
    \partial_t \Ent{h}{g} = \int \left(\partial_t h\right) \ln(\frac{h}{g}) + h \partial_t\left(\ln\left(\frac{h}{g}\right)\right).
 \end{equation*}

Replacing both $h_t$ and $g_t$ by their dynamics and integrating by parts  shows

\begin{equation*}
    \begin{aligned}
        \partial_t \Ent{h}{g} &= \int \left(\lap h + \grad(hB) \right)\ln\left(\frac{h}{g}\right) + h\partial_t\left(\grad \ln(h) - \grad\ln(g)\right) \\
        &= \int  \left(\lap h + \grad(hB) \right)\ln\left(\frac{h}{g}\right) + h\left( \frac{\lap h + \grad (hB)}{h} - \frac{\lap g + \grad(gA)}{g} \right) \\
        &= \int  \left(\lap h + \grad(hB) \right)\ln\left(\frac{h}{g}\right) + \lap h + \grad(hB) -\frac{h}{g}\grad \cdot \left(\grad g + gA \right) \\
        &= \int \grad \cdot (\grad h + hB)(\ln\left(\frac{h}{g}\right) + 1) - \int \frac{h}{g}\grad \cdot (\grad g + gA) \\
        &= -\int \left(\grad h + hB \right) \cdot \grad \ln\left(\frac{h}{g}\right) + \int \grad\left(\frac{h}{g}\right)\cdot(\grad g + gA).
    \end{aligned}
\end{equation*}

Using $\grad h = h \grad \ln h$ and adding manually $h \grad \ln(g)$ entails

\begin{equation*}
\begin{aligned}
    \partial_t \Ent{h}{g} &= -\int h \grad \ln(h)\cdot\grad\ln\left(\frac{h}{g}\right) - \int hB\grad\ln\left(\frac{h}{g}\right) - \grad(\frac{h}{g})\cdot\left( \grad g + gA \right)\\
    &= -\int \left(h \grad \ln(h) - h \grad \ln(g)\right)\cdot \grad \ln\left(\frac{h}{g}\right) \\ & \quad - \int (h \grad \ln(g) + hB)\cdot \grad \ln\left(\frac{h}{g}\right) + \grad\left(\frac{h}{g}\right)\cdot(\grad g + gA) \\
    &= -\mathcal{I}(h|g) - \int \grad \ln\left(\frac{h}{g}\right)\cdot\left(h\grad \ln(g) + hB\right) - \grad \left(\frac{h}{g}\right)\cdot\left(\grad g + gA\right).
\end{aligned}
\end{equation*}
Now, we use $\grad\ln\left(\frac{h}{g}\right)=\frac{g\grad\left(\frac{h}{g}\right)}{h}$ and $g\grad\ln(g) = \grad g$:

\begin{equation*}
    \begin{aligned}
        \partial_t \Ent{h}{g}&=-\mathcal{I}(h|g)-\int \grad\left(\frac{h}{g}\right)\left(g\grad\ln(g) + gB - \grad g - gA\right) \\
        &= -\mathcal{I}(h|g) - \int g\grad\left(\frac{h}{g}\right)(B-A) \\
        &= - \mathcal{I}(h|g) - \int h\grad\ln\left(\frac{h}{g}\right)(B-A).
    \end{aligned}
\end{equation*}
\end{eproof}


Applying Proposition \ref{dynamiquentropie} to the dynamics of $u_t^{(N)}$ and $u_t^{\otimes N}$ as given in Proposition \ref{dynamiques} yields

\begin{equation}\label{dynamique1}
    H'(t) = -\int \Big|\frac{\grad u_t^{(N)}}{u_t^{(N)}} - \frac{\grad u_t^{\otimes N}}{u_t^{\otimes N}}\Big|^2u_t^{(N)} - \int \left( \frac{\grad u_t^{(N)}}{u_t^{(N)}} - \frac{\grad u_t^{\otimes N}}{u_t^{\otimes N}}\right)\cdot[\grad \psi^{(N)} - \grad \psi_t ]u_t^{(N)}    
\end{equation}

where $\psi^{(N)}$ was defined in Proposition \ref{dynamiques} and $\psi_t : \R^{dN} \longrightarrow \R$ is defined for $(x_1\cdots,x_N) \in \R^{dN}$ by 
\begin{equation*}
    \psi_t(x_1,\cdots,x_N) = \sum_{i=1}^N V(x_i) +  W \ast u_t (x_i).
\end{equation*}
In particular 
 $(\grad \psi_t)_i = \grad V(x_i) + \grad W \ast u_t(x_i)  $ and $(\grad \psi^{(N)})_i = \grad V(x_i) + \frac{1}{N}\sum_{j=1}^N \grad W (x_i - x_j)$.

This brings us to a differential inequality on $H$, as in \cite{Yau91}, \cite{OVa91} and \cite{malrieu}.
 \begin{eprop}\label{termeerreur}
     We have the following differential inequality 

     \begin{equation}\label{INEQDIFF}
         H'(t) \leq -\frac{1}{2\CLSI}H(t) + \frac{1}{2}\mathcal{E}_N(t)
     \end{equation}

     where the error term $\mathcal{E}_N(t)$ is defined as 

     \begin{equation*}
       \mathcal{E}_N(t) :=  \int\abs{\grad \psi^{(N)} - \grad \psi_t}^2u_t^{(N)} = \sum_{i=1}^N\E{\abs{\grad W \ast u_t(X^{i,N}) - \frac{1}{N}\sum_{j=1}^N \grad W(X^{i,N} - X^{j,N}) }^2}
     \end{equation*}

 \end{eprop}

 \begin{eproof}
     We begin by using $a\cdot b \leq \frac{1}{2}\left(a^2 + b^2\right)$ on the right-hand side of \eqref{dynamique1} to prove

     \begin{equation*}
     \begin{aligned}
     H'(t) &\leq -\frac{1}{2}\int \Big|\frac{\grad u_t^{(N)}}{u_t^{(N)}} - \frac{\grad u_t^{\otimes N}}{u_t^{\otimes N}}\Big|^2u_t^{(N)} + \frac{1}{2}\int\abs{\grad \psi^{(N)} - \grad \psi_t}^2u_t^{(N)}\\
     &= -\frac{1}{2}\mathcal{I}\left(u_t^{(N)}|u_t^{\otimes N}\right) + \frac{1}{2}\mathcal{E}_N(t).
     \end{aligned}
     \end{equation*}

Moreover, we have a uniform-in-time log-Sobolev inequality for $u_t$ with constant $\CLSI$ thanks to Proposition \ref{propaglsilimit} and Remark \ref{remark30}. Then, by tensorization we have a log-Sobolev inequality for $u_t^{\otimes N}$ with the same constant (see \cite[Theo. 3.22]{logsob} or \cite[Prop. 5.2.7]{BGL2014})
Therefore 

\begin{equation*}
    \mathcal{I}\left(u_t^{(N)}|u_t^{\otimes N}\right) \geq \frac{1}{\CLSI}H(t),
\end{equation*}
which concludes.
 \end{eproof}

To turn \eqref{INEQDIFF} into a control on $H(t)$ using Grönwall's lemma, we need to find an upper bound on $\mathcal{E}_N(t)$. All the proofs follow the same scheme, however the regularity of $W$ will lead to different control. We therefore give a lot of detail for the Lipschitz case (Theorem \ref{cas lipschitz}), to which the other proofs will refer.

\subsection{Proof of Theorem \ref{cas lipschitz}}
\begin{elem}\label{diffgradlip}
    Under $\Hypp$, assume that $\grad W$ is Lipschitz. Then,  for all $t \geq 0$
    \begin{equation*}
            \mathcal{E}_N(t) \leq \Constinit\left(e^{-2\beta t}\mathbf{W}_2^2\left(u_0^{(N)},u_0^{\otimes N}\right) + 1\right).
    \end{equation*}
\end{elem}
\bigskip
\begin{eproof}[Lemma \ref{diffgradlip}]
In what follows, we omit the subscript $"t"$ on the particles.
Let $L_W$ be a Lipschitz constant for $\grad W$. Let $1 \leq i \leq N$. Then 
\begin{equation}\label{ineqtriangle}
\begin{aligned}
    \abs{\grad W \ast u_t(X^{i,N}) - \frac{1}{N}\sum_{j=1}^N   \grad W (X^{i,N} - X^{j,N})}^2 &\leq
           3 \abs{\grad W \ast u_t (X^{i,N}) - \grad W\ast u_t (\bar{X^i})}^2 \\  &\quad + 3\abs{\grad W \ast u_t(\bar{X^i}) - \frac{1}{N}\sum_{j=1}^N   \grad W (\bar{X^i} - \bar{X^j}) }^2  \\ &\quad + 3\abs{\frac{1}{N}\sum_{j=1}^N \grad W (\bar{X^i} - \bar{X^j}) - \grad W (X^{i,N} - X^{j,N}) }^2 \\
           &:=  \mathbf{I} + \mathbf{II} + \mathbf{III}
\end{aligned}
        \end{equation}
We begin with $\mathbf{I}$. The Lipschitz regularity of $\grad W$ yields
\begin{equation*}
 \abs{\grad W \ast u_t (X^{i,N}) - \grad W\ast u_t (\bar{X^i})}^2 \leq L_W^2\abs{X^{i,N}-\bar{X^i}}^2.
\end{equation*}

Summing on $i$ entails
\begin{equation}\label{bigstar}
 \sum_{i=1}^N\abs{\grad W \ast u_t (X^{i,N})- \grad W\ast u_t (\bar{X^i})}^2  \leq L_W^2\sum_{i=1}^N\abs{X^{i,N}-\bar{X^i}}^2.
\end{equation}

 We continue with $\mathbf{III}$. Using  Jensen's inequality, the fact that $(a+b)^2 \leq 2(a^2 + b^2)$ and the Lipschitz regularity of $\grad W$ yields
\begin{equation*}
    \begin{aligned}
     \abs{\frac{1}{N}\sum_{j=1}^N \grad W (\bar{X}^i - \bar{X^j}) - \grad W (X^{i,N} - X^{j,N}) }^2 &\leq \frac{L_W^2}{N} \sum_{j=1}^N \abs{\bar{X^i} - \bar{X^j} + X^{j,N}- X^{i,N}}^2\\ 
     &\leq \frac{2L_W^2}{N}\sum_{j=1}^N \left(\abs{\bar{X^i} - X^{i,N}}^2 + \abs{\bar{X^j} - X^{j,N}}^2\right) \\
     & = \frac{2L_W^2}{N} \sum_{j=1}^N \abs{\bar{X^j} - X^{j,N}}^2 + 2L_W^2 \abs{\bar{X^i} - X^{i,N}}^2.
    \end{aligned}
\end{equation*}

Summing on $i$ yields

\begin{equation}\label{RHS}
     \sum_{i=1}^N\abs{\frac{1}{N}\sum_{j=1}^N \grad W (\bar{X}^i - \bar{X^j}) - \grad W (X^{i,N} - X^{j,N}) }^2 \leq 4L_W^2\sum_{i=1}^N \abs{\bar{X^i} - X^{i,N}}^2.
\end{equation}

We have already treated the term $\mathbf{II}$ in the proof of Theorem \ref{thcouplage} with a law of large numbers argument (Step 2.1), where we found that \begin{equation}\label{***}
\sum_{i=1}^N \E{\abs{\grad W \ast u_t(\bar{X^i}) -\frac{1}{N}\sum_{j=1}^N  \grad W (\bar{X^i} - \bar{X^j})}^2} \leq \Constinit.
\end{equation}
Injecting  \eqref{bigstar},\eqref{RHS}, \eqref{***} in \eqref{ineqtriangle} and using Theorem \ref{thcouplage} entails

\begin{equation*}
    \mathcal{E}_N(t) \leq \Constinit\left( 1 + \sum_{i=1}^N \E{\abs{\bar{X^i} - X^{i,N}}^2}\right) 
    \leq \Constinit\left(1 +e^{-2\beta t}\sum_{i=1}^N  \E{\abs{X^{i,N}_0 - \bar{X^i_0}}^2} \right).
    \end{equation*}

Taking the infimum on the initial conditions yields

        \begin{equation*}
            \mathcal{E}_N(t) \leq \Constinit\left(e^{-2\beta t}\mathbf{W}_2^2\left((u_0^{(N)},u_0^{\otimes N}\right) + 1\right).
    \end{equation*}
    \end{eproof}
\begin{eproof}[Theorem \ref{cas lipschitz}]

Using Lemma \ref{diffgradlip} in \eqref{INEQDIFF} and using Grönwall's lemma shows that 
    \begin{equation*}
               \frac{ \Ent{u_t^{(N)}}{u_t^{\otimes N}}}{N} \leq  e^{-\frac{1}{2\CLSI}t}\Ent{u_0^{(N)}}{u_0^{\otimes N}} + \frac{\Constinit}{N}\left(e^{-\gamma t}\mathbf{W}_2^2\left(u_0^{(N)},u_0^{\otimes N}\right)  + 1\right)
    \end{equation*}

Since we have a log-Sobolev inequality with constant $\CLSI$ for $u_0$ and by tensorization with the same constant for $u_0^{\otimes N}$, we have a Talagrand's inequality with the same constant for $u_0^{\otimes N}$ (see \cite[Theorem 1]{OV}). 
This shows

    \begin{equation*}
               \frac{ \Ent{u_t^{(N)}}{u_t^{\otimes N}}}{N} \leq    \frac{\Constinit}{N}\left(e^{-\gamma t}\Ent{u_0^{(N)}}{u_0^{\otimes N}}  + 1\right),
    \end{equation*}

which concludes the proof.

\end{eproof}
\begin{erem}\label{remarquemoments}
It is natural to try and extend this method to the case where $\grad W$ is not Lipschitz but with polynomial growth.

Using our moments control from Proposition \ref{controlemomentsparticule} and \ref{particulelimiteborne}, the same procedure applies to control \textbf{II}. However, when trying to control \textbf{III}, the situation is quite different. For instance, the same procedure on the term~$\mathbf{I}$ gives the following 

\begin{equation*}
\E{\abs{ \grad W \ast u_t(X^{i,N})
        - \grad W \ast u_t(\bar{X^i}) }^2}
\\\leq M^2  
\E{\abs{X^{i,N} - \bar{X^i}}^2
\left(1 + \abs{X^{i,N}}^{2m} 
        + \abs{\bar{X^i}}^{2m}\right)}.
\end{equation*}
The issue here is the emergence of terms such as $\E{\abs{X^{i,N} - \bar{X^i}}^2\abs{X^{i,N}}^{2m}}$, on which we have no control \textit{a priori}. 
We resort to the Hölder inequality to prove, with $0 < a < 2$ 
\begin{equation*}
    \E{\abs{X^{i,N} - \bar{X^i}}^2\abs{X^{i,N}}^{2m}} \leq  \left(\E{\abs{X^{i,N} - \bar{X^i}}^2}\right)^{\frac{a}{2}}\left(\E{\abs{X^{i,N} - \bar{X^i}}^2\abs{X^{i,N}}^{\frac{4m}{2-a}}}\right)^{\frac{2-a}{2}}.
\end{equation*}
The last factor can be bounded using the moment estimates, provided the initial conditions admit sufficently high moments.
However, the fact that we cannot take $a = 2$ implies that this control might not be enough to reach a rate in $1/N$. For instance, in the \textit{i.i.d} and log-Sobolev setting for the particle system, we have only been able to derive a control  of $\Ent{u_t^{(N)}}{u_t^{\otimes N}}/N$ in $1/N^\alpha$ for all $0 < \alpha (= \frac{a}{2}) < 1$. 
This control in $\frac{1}{N^\beta}$ can be improved by a Markov-like argument, to a control in $\ln(N)^m/N$ where $m$ is the degree of the polynomial growth of $\grad W$, as long as we have a uniform in $i,t,N$ Gaussian moment on $\abs{X_t^i}$ (that is, in particular, the case when we have a log-Sobolev inequality with a constant that does not depend on $N$). 
\end{erem}

\subsection{Proof of Theorem \ref{th4}}
\begin{elem}\label{prop37}
    Under $\Hypp$, assume that for all $p > 0$, $\displaystyle \sup_{1 \leq i \leq N, N \in \N} \int \abs{x_i}^{2p} \mathrm{d}u_0^{(N)} < + \infty$ and that $u_0 \in \mathcal{P}_{2p}(\R^d)$ for all $p \in \N^*$. 
    Then, for all $t \geq 0$, $0 < \alpha < 1$:
        \begin{equation*}
    \mathcal{E}_N(t) \leq \Constinit(\alpha)N^{1-\alpha}\left(1 +  e^{-2\beta \alpha t}\mathbf{W}_2^{2\alpha}\left((u_0^{(N)},u_0^{\otimes N}\right)\right).
\end{equation*}
\end{elem}
\begin{erem}
    If we only suppose $u_0 \in \mathcal{P}_{2mp}(\R^d)$ and $\displaystyle \sup_{1 \leq i \leq N, N \in \N} \int \abs{x_i}\mathrm{d}u_0^{(N)} < + \infty$ for some $p \in \N^*$, we get the result for $0 < \alpha < 1 - \frac{m}{p-1}$.
\end{erem}

\begin{eproof}[Lemma \ref{prop37}]
    We start anew from \eqref{ineqtriangle} and the control of terms \textbf{I}, $\mathbf{II}, \mathbf{III}$. 
    As we have initial moments of order $2m$, we can control the term $\mathbf{II} $, that is \begin{equation*}
      \mathbf{II} = \sum_{i=1}^N\E{\abs{\grad W \ast u_t (\bar{X^i}) -\frac{1}{N}\sum_{j=1}^N  \grad W (\bar{X^i} - \bar{X^j})}^2} \leq \Constinit
    \end{equation*} 
    with a law of large numbers argument as in the proof of Theorem \ref{thcouplage} thanks to Proposition \ref{particulelimiteborne}.

The issue here is the control of $\mathbf{I}$ and $\mathbf{III} $. For instance, we begin by $\mathbf{III}$: \begin{equation*}
     \textbf{III} = \sum_{i=1}^N\E{\abs{\frac{1}{N}\sum_{j=1}^N \grad W (\bar{X^i} - \bar{X^j}) - \grad W (X^{i,N} - X^{j,N}) }^2}.
\end{equation*}

Let $0 < a < 2$. Using Jensen's inequality on the sum on $j$, polynomial growth and Hölder inequality, we obtain

\begin{equation}\label{holder}
\begin{gathered}
         \mathbf{III} \leq \frac{\Const}{N}\sum_{i,j=1}^N\E{\abs{\bar{X^i}  - \bar{X^j} - X^{i,N} + X^{j,N} }^2\left(1 +  \abs{\bar{X^i} - \bar{X^j}}^{2m} + \abs{X^{i,N}-X^{j,N}}^{2m} \right)} \\
         \leq \frac{\Const}{N}\sum_{i,j=1}^N \E{\abs{X^{j,N} - \bar{X^j}}^2} +\frac{\Const}{N} \sum_{i,j=1}^N  \left(\E{\abs{X^{j,N} - \bar{X^j}}^2}\right)^{\frac{a}{2}}\left(\E{\abs{X^{j,N} - \bar{X^j}}^2\abs{\bar{X^i} - \bar{X^j}}^{\frac{4m}{2-a}}}\right)^{\frac{2-a}{2}}  \\ +   \frac{\Const}{N}\sum_{i=1}^N\sum_{j=1}^N \left(\E{\abs{X^{j,N} - \bar{X^j}}^2}\right)^{\frac{a}{2}}\left(\E{\abs{X^{j,N} - \bar{X^j}}^2\abs{X^{i,N} - X^{j,N}}^{\frac{4m}{2-a}}}\right)^{\frac{2-a}{2}}.
\end{gathered}
\end{equation}

Since we have moments of all orders, using Hölder inequalities and our moment control for both the limit particle and the particle system from Proposition \ref{controlemomentsparticule} and \ref{particulelimiteborne} shows us that this is all bounded by 
\begin{equation}\label{borne1}
    \Const \sum_{j=1}^N \E{\abs{X^{j,N} - \bar{X^j}}^2} + \Constinit(\alpha)\sum_{j=1}^N \left(\E{\abs{X^{j,N} - \bar{X^j}}^2}\right)^\alpha
\end{equation}
for all $0 < \alpha < 1$ (we took $\alpha = \frac{a}{2}$).
We can therefore deteriorate the constant $\Constinit(\alpha)$ with moments of order 2 with the following remark:

\begin{equation*}
    \begin{aligned}
    \E{\abs{X^{j,N} - \bar{X^j}}^2} &= \left(\E{\abs{X^{j,N}-\bar{X^j}}^2}\right)^\alpha\left(\E{\abs{X^{j,N}-\bar{X^j}}^2}\right)^{1- \alpha} \\ &\leq \left(\E{\abs{X^{j,N}-\bar{X^j}}^2}\right)^\alpha\left(2^{1- \alpha}\left(\E{\abs{X^{j,N}}^2} + \E{\abs{\bar{X^j}}^2}\right)^{1-\alpha}\right) \\&
    \leq \Constinit(\alpha)\left(\E{\abs{X^{j,N}-\bar{X^j}}^2}\right)^\alpha.
    \end{aligned}
\end{equation*}

Here, we used the fact that $0 < \alpha < 1$ which ensures that $x \mapsto x^{1-\alpha}$ is increasing, the inequality $(a+b)^2 \leq 2(a^2 + b^2)$ and in the last inequality we used the propagation of moments through Proposition \ref{controlemomentsparticule} and \ref{particulelimiteborne}.
Moreover, $x \mapsto x^\alpha$ is concave since $0 < \alpha < 1$. Hence, Jensen's inequality yields, for all positive $(x_i)$,
\begin{equation*}
    \sum_{j=1}^N x_j^\alpha \leq N^{1-\alpha} \left(\sum_{j=1}^N x_j\right)^\alpha.
\end{equation*}

Injecting our last two remarks in \eqref{borne1}, we obtain 

\begin{equation*}
\begin{aligned}    
    \sum_{i=1}^N\E{\abs{\frac{1}{N}\sum_{j=1}^N \grad W (\bar{X^i} - \bar{X^j}) - \grad W (X^{i,N} - X^{j,N}) }^2} &\leq \Constinit(\alpha)N^{1-\alpha}\left(\sum_{j=1}^N \E{\abs{X^{j,N} - \bar{X^j}}^2}\right)^\alpha .
    \end{aligned}
\end{equation*}

We have now bounded $\mathbf{III}$ and $\mathbf{II}$. The term $\mathbf{I}$ can be controlled just like $\mathbf{III}$ by a similar term.
Combining our controls of \textbf{I}, \textbf{II} and \textbf{III} in \eqref{ineqtriangle} shows that 
\begin{equation*}
    \mathcal{E}_N(t) \leq \Constinit\left(1 + N^{1-\alpha}\left(\sum_{j=1}^N\E{\abs{X^{j,N} - \bar{X^j}}^2}\right)^\alpha\right).
\end{equation*}

Using Theorem \ref{coupling} shows, for all $0 < \alpha < 1$

\begin{equation*}
    \mathcal{E}_N(t) \leq \Constinit(\alpha)N^{1-\alpha}\left(1 +  e^{-2\beta \alpha t}\left(\sum_{i=1}^N  \E{\abs{X^{i,N}_0 - \bar{X^i_0}}^2}\right)^\alpha\right).
\end{equation*}

Taking the infimum on coupling of initial conditions concludes.
\bigskip
\end{eproof}
\begin{eproof}[Theorem \ref{th4}]
Plugging Lemma \ref{prop37} in \eqref{INEQDIFF}, yields

\begin{equation}\label{retour alphacritique}
             H'(t) \leq \frac{-1}{2\CLSI}H(t) + \Constinit(\alpha)N^{1-\alpha}\left(1 +  e^{-2\beta \alpha t}\mathbf{W}_2\left(u_0^{(N)},u_0^{\otimes N}\right)^{2\alpha}\right).
\end{equation}

Applying Grönwall's lemma for $\alpha \neq \frac{1}{4\beta \CLSI}= \alpha_C$ entails

\begin{equation*}
    H(t) \leq e^{-\frac{1}{2\CLSI}t}H(0) + \Constinit(\alpha) N^{1-\alpha}\left(1 + \mathbf{W}_2\left(u_0^{(N)},u_0^{\otimes N}\right)^{2\alpha}\frac{e^{-2\beta \alpha t}-e^{\frac{-1}{2C}t}}{\frac{1}{2\CLSI}-2\beta \alpha} \right).
    \end{equation*}

that is 

\begin{equation*}
    \frac{H(t)}{N} \leq e^{-\frac{1}{2\CLSI}t}\frac{H(0)}{N} +  \frac{\Constinit(\alpha)}{N^\alpha}\left(1 + e^{-\gamma t}\mathbf{W}_2^{2\alpha}(u_0^{N},u_0^{\otimes N})\right)
\end{equation*}

Now, we can use the fact that $x^{\alpha} \leq 1 + x$ for $\alpha \leq 1$ with $x = \mathbf{W}_2^2\left(u_0^{(N)},u_0^{\otimes N}\right)$ and a Talagrand's inequality for $u_0^{\otimes N}$ implied by the log-Sobolev inequality on $u_0$ and hence $u_0^{\otimes N}$ to prove that 

\begin{equation*}
    \frac{H(t)}{N} \leq \frac{\Constinit(\alpha)}{N^\alpha}\left(1 + e^{-\gamma t}H(0)\right),
\end{equation*}
which concludes.
The result remains true for $\alpha = \alpha_C$. One just has to use in \eqref{retour alphacritique} the fact that for all $t \geq 0$, $\alpha < \alpha_C < 1, e^{-2\beta\alpha_C t} \leq e^{-2\beta \alpha t}$.  Then, taking $\Constinit(\alpha_C) = \Constinit(\alpha)$ concludes.
\end{eproof}
\begin{erem}
    We cannot take the limit $\alpha \longrightarrow 1$ as $\Constinit(\alpha)$ might  blow up as $\alpha$ goes to $1$: we do not immediately have a convergence rate of $\frac{1}{N}$ as in the Lipschitz case.
\end{erem}
\subsection{Proof of Theorem \ref{log}}

\begin{elem}\label{resultatfinal}
  Under $\Hypp$, assume that $u_0^{(N)}$ satisfies a log-Sobolev inequality with constant independent of $N$ and that $\displaystyle \sup_{1 \leq i \leq N,N \in \N} \int\abs{x_i}^{2m^2}\mathrm{d}u_0^{(N)} < + \infty$. In addition, assume that $u_0$ satisfies a log-Sobolev inequality with constant $c_0$.
    Then, for all $t \geq 0$:
        \begin{equation*}
           \mathcal{E}_N(t) \leq \Constinit\ln(N)^m\left(1 + e^{-2\beta t}\mathbf{W}_2^2\left(u_0^{(N)},u_0^{\otimes N}\right) \right).
    \end{equation*}
    \end{elem}
\bigskip

\begin{erem}
    In Theorem \ref{log}, we need a logarithmic Sobolev inequality for both the limit system and the particle system, which contrasts with our previous results which needed a logarithmic Sobolev inequality only on the limit system. Here, the logarithmic Sobolev inequality on the particle system is only used to propagate sub-Gaussian integrability.
\end{erem}
\begin{eproof}[Lemma \ref{resultatfinal}]
First, note that since $u_0^{(N)}$ satisfies a logarithmic Sobolev inequality, according to Proposition \ref{proplsiparticle}, $u_t^{(N)}$ satisfies a logarithmic Sobolev inequality with constant $\CLSI$.

It is well known that a logarithmic Sobolev inequality implies some exponential integrability. It is less well-known that this exponential moment can be bounded by a quantity that only depends on the second moment of the measure. Indeed, thanks to the proof of \cite[Proposition 5.4.1]{BGL2014}, notably the last equation of the proof, there exists $\alpha > 0$ such that \begin{equation}\label{ControleLSI}
    \sup_{t,i,N} \E{e^{\alpha \abs{X_t^{i,N}}^2}} \leq \Const(\CLSI, \E{\abs{X_t^{i,N}}^2}) \leq \Constinit.
\end{equation}
Where the last inequality is due to the fact that $ \displaystyle \sup_{t,i,N} \E{ \abs{X_t^{i,N}}^2} < + \infty$  thanks to Proposition \ref{controlemomentsparticule}.
Note that the constant $\Constinit$ also depends on $\CLSI$.
The same reasoning with the limit system shows that both our systems are exponentially integrable, which entails that they both have moments of any order, which are uniformly bounded in time, $i$ and $N$ thanks to Proposition \ref{controlemomentsparticule} and \ref{particulelimiteborne}.

Starting anew from \eqref{holder}, we have to control 

    \begin{equation*}
        \sum_{1 \leq i,j \leq N}\E{\abs{X_t^{j,N} - \bar{X}_t^j - X_t^{i,N} + \bar{X}_t^i}^2\left(1 + \abs{\bar{X}_t^i - \bar{X}_t^j}^{2m} + \abs{X_t^{i,N} - X_t^{j,N}}^{2m}\right)}.
        \end{equation*}

Using $(a+b)^2 \leq 2(a^2+b^2),$ moment control through Proposition \ref{controlemomentsparticule} and \ref{particulelimiteborne} shows that it is bounded by 

\begin{equation*}
   \Constinit \sum_{i=1}^N \E{\abs{X_t^{i,N}-\bar{X}_t^i}^2\left(1 + \abs{\bar{X}_t^i}^{2m} + \abs{X_t^{i,N}}^{2m}\right)}.
\end{equation*}

It remains to study terms such as 
\begin{equation}\label{les2termes}
 \sum_{i=1}^N\E{\abs{X_t^{i,N}-\bar{X}_t^i}^2\abs{\bar{X}_t^i}^{2m}} \quad \text{and} \quad    \sum_{i=1}^N \E{\abs{X_t^{i,N}-\bar{X}_t^i}^2\abs{X_t^{i,N}}^{2m}}.
\end{equation}

We will only study the second one, the procedure being the same for the first one. Let $R > 0$.

\begin{equation}\label{markov1}
\begin{aligned}
          \sum_{i=1}^N \E{\abs{X_t^{i,N}-\bar{X}_t^i}^2\abs{X_t^{i,N}}^{2m}} &=       \sum_{i=1}^N \E{\abs{X_t^{i,N}-\bar{X}_t^i}^2\abs{X_t^{i,N}}^{2m}\mathbf{1}_{\abs{X_t^{i,N}} \leq R}} \\ &\quad+   \sum_{i=1}^N\E{\abs{X_t^{i,N}-\bar{X}_t^i}^2\abs{X_t^{i,N}}^{2m}\mathbf{1}_{\abs{X_t^{i,N}}> R}}.   
\end{aligned}
\end{equation}

We begin with the first term in \eqref{markov1}.

\begin{equation*}
    \sum_{i=1}^N \E{\abs{X_t^{i,N}-\bar{X}_t^i}^2\abs{X_t^{i,N}}^{2m}\mathbf{1}_{\abs{X_t^{i,N}} \leq R}} \leq R^{2m}\sum_{i=1}^N \E{\abs{X_t^{i,N} - \bar{X}_t^i}^2}.
\end{equation*}

Now, we study the second term in \eqref{markov1}.

Let $0 < 2 \varepsilon < \alpha $. For $R$ large enough, $y \mapsto y^{2m}e^{-\varepsilon y^2}$ is decreasing on $[R,+\infty)$. We choose such an $R$ and therefore, for all $1 \leq i \leq N$,
\begin{equation*}
    \E{\abs{X_t^{i,N}-\bar{X}_t^i}^2\abs{X_t^{i,N}}^{2m}\mathbf{1}_{\abs{X_t^{i,N}}> R}} \leq R^{2m}e^{-\varepsilon R^2} \E{\abs{X_t^{i,N}- \bar{X}_t^i}^2e^{\varepsilon \abs{X_t^{i,N}}^2}}.
\end{equation*}

Injecting this in the second term of \eqref{markov1} and using Cauchy-Schwarz's inequality shows that the second term of \eqref{markov1} is bounded by

\begin{equation}\label{markov12}
R^{2m}e^{-\varepsilon R^2}\left(\sum_{i=1}^N \E{\abs{X_t^{i,N} - \bar{X}_t^i}^4}\right)^{\frac{1}{2}}\left(\sum_{i=1}^N \E{e^{2\varepsilon \abs{X_t^{i,N}}^2}}\right)^{\frac{1}{2}}.
\end{equation}

Using propagation of moments with Proposition \ref{controlemomentsparticule} and \ref{particulelimiteborne}  on the first term of the product and sub-Gaussian integrability through \eqref{ControleLSI} on the second term shows that \eqref{markov12} is bounded by 

\begin{equation*}
    \Constinit R^{2m}e^{-R^2}N.
\end{equation*}

Therefore, \eqref{markov1} is bounded by 

\begin{equation*}
\Constinit R^{2m}\left( \sum_{i=1}^N \left(\E{\abs{X_t^{i,N}- \bar{X}_t^i}^2} \right) + Ne^{-\varepsilon R^2}\right).
\end{equation*}

Choosing $R^2 = \frac{\ln(N)}{\varepsilon}$ yields

\begin{equation*}
     \sum_{i=1}^N \E{\abs{X_t^{i,N}-\bar{X}_t^i}^2\abs{X_t^{i,N}}^{2m}} \leq \Constinit \ln(N)^m\left(1 + \sum_{i=1}^N \E{\abs{X_t^{i,N} - \bar{X}_t^i}^2}\right).
\end{equation*}

The first term of \eqref{les2termes} is identically handled to prove

\begin{equation*}
     \sum_{i=1}^N \E{\abs{X_t^{i,N}-\bar{X}_t^i}^2\abs{\bar{X^i_t}}^{2m}} \leq \Constinit \ln(N)^m\left(1 + \sum_{i=1}^N \E{\abs{X_t^{i,N} - \bar{X}_t^i}^2}\right).
\end{equation*}
Therefore, using Theorem \ref{thcouplage} yields

    \begin{equation*}
           \mathcal{E}_N(t) \leq \Constinit\ln(N)^m\left(e^{-2\beta t}\sum_{i=1}^N  \E{\abs{X^{i,N}_0 - \bar{X^i_0}}^2} + 1\right).
    \end{equation*}
    \end{eproof}
\begin{eproof}[Theorem \ref{log}]
As in the proof of Theorem \ref{th4}, using Lemma \ref{resultatfinal} and Grönwall's lemma shows that
\begin{equation*}
    \frac{\Ent{u_t^{(N)}}{u_t^{\otimes N}}}{N} \leq e^{-\frac{1}{2\CLSI}t}  \frac{\Ent{u_0^{(N)}}{u_0^{\otimes N}}}{N} + \Constinit\frac{\ln(N)^m}{N}\left(1 + e^{-\gamma t}\mathbf{W}_2^2\left(u_0^{(N)},u_0^{\otimes N}\right)\right).
    \end{equation*}

using a Talagrand's inequality as in the proof of Theorem \ref{cas lipschitz} concludes.
\end{eproof}


\subsection{Proof of Theorem \ref{thd=1}}
\begin{elem}\label{lemmed1}
 Let $p \geq 2, p \in \N$. Under $\Hypp$, assume that $\displaystyle\sup_{1 \leq i \leq N, N \in \N} \int \abs{x_i}^{2m^2p} \mathrm{d}u_0^{(N)}<+\infty$, $u_0 \in \mathcal{P}_{2m^2p}(\R^d)$.
 If either $W = \lambda_W \frac{\abs{x}^2}{2}$, or $ d = 1$, then, for all $t \geq 0$:
    \begin{equation*}
          \mathcal{E}_N(t) \leq \Constinit\left(1 + e^{-2\beta t} N^{1-\frac{1}{p}}\mathbf{W}_{2p}^2\left((u_0^{(N)},u_0^{\otimes N}\right) \right).
    \end{equation*}

  \end{elem}
\begin{eproof}[Lemma \ref{lemmed1}]
    We start anew from \eqref{holder}. As said in the proof of Theorem \ref{log}, we have to bound terms of the form 

    \begin{equation*}
            \sum_{i=1}^N \E{\abs{X_t^{i,N}-\bar{X}_t^i}^2\abs{X_t^{i,N}}^{2m}} \quad \text{and} \quad \E{\abs{X_t^{i,N}-\bar{X}_t^i}^2\abs{\bar{X}_t^i}^{2m}}.
    \end{equation*}

We only treat the first term. The difference with \eqref{holder} is that we can use Hölder's inequality with an exponent $p\geq2$. Indeed, Hölder's inequality yields

\begin{equation*}
        \sum_{i=1}^N \E{\abs{X_t^{i,N}-\bar{X}_t^i}^2\abs{X_t^{i,N}}^{2m}} \leq \left(\sum_{i=1}^N \E{\abs{X_t^{i,N}- \bar{X}_t^i}^{2p}}\right)^{\frac{1}{p}}\left(\sum_{i=1}^N \E{\abs{X_t^{i,N}}^{\frac{2mp}{p-1}}}\right)^{1-\frac{1}{p}}.
    \end{equation*}

Using Theorem \ref{thcouplage} through \eqref{wass} and propagation of moments through Proposition \ref{controlemomentsparticule} with $k = m\frac{p}{p-1}$ yields 

\begin{equation*}
\begin{aligned}
\sum_{i=1}^N \E{\abs{X_t^{i,N}-\bar{X}_t^i}^2\abs{X_t^{i,N}}^{2m}} &\leq \Constinit N^{1-\frac{1}{p}} \left(e^{-2\beta t }\left(\sum_{i=1}^N \E{\abs{X_0^{i,N}- \bar{X}_0^i}^{2p}}\right)^{\frac{1}{p}} + \frac{1}{N^{1-\frac{1}{p}}}\right)  \\
        & = \Constinit \left(N^{1-\frac{1}{p}}e^{-2\beta t}\left(\sum_{i=1}^N \E{\abs{X_0^{i,N}- \bar{X}_0^i}^{2p}}\right)^{\frac{1}{p}} + 1\right).
\end{aligned}
\end{equation*}

Collecting all the terms shows that
    \begin{equation*}
         \mathcal{E}_N(t) \leq   \Constinit\left(e^{-2\beta t}\left(\sum_{i=1}^N  \E{\abs{X^{i,N}_0 - \bar{X^i_0}}^2} + N^{1-\frac{1}{p}}\left(\sum_{i=1}^N \E{\abs{X_0^{i,N}- \bar{X}_0^i}^{2p}}\right)^{\frac{1}{p}}\right)+ 1 \right).
    \end{equation*}

 Using Hölder's inequality twice, once on $\E{\abs{X_0^{i,N}-\bar{X}_0^i}^2}$ and the second time on the sum yields 
    \begin{equation*}
          \mathcal{E}_N(t) \leq \Constinit\left(e^{-2\beta t} N^{1-\frac{1}{p}}\left(\sum_{i=1}^N \E{\abs{X_0^{i,N}- \bar{X}_0^i}^{2p}}\right)^{\frac{1}{p}}+ 1 \right).
    \end{equation*} \end{eproof}
    \begin{eproof}[Theorem \ref{thd=1}]

    As in the proof of Theorem \ref{th4}, using Grönwall's lemma yields

    \begin{equation*}
        \frac{\Ent{u_t^{(N)}}{u_t^{\otimes N}}}{N} \leq e^{-\frac{1}{2\CLSI}t}  \frac{\Ent{u_0^{(N)}}{u_0^{\otimes N}}}{N} + \frac{\Constinit}{N}\left( 1 + e^{-\gamma t}N^{1-\frac{1}{p}}\left(\sum_{i=1}^N \E{\abs{X_0^{i,N}- \bar{X}_0^i}^{2p}}\right)^{\frac{1}{p}}\right).
    \end{equation*}
Taking the infimum on the initial conditions yields
\begin{equation*}
        \frac{\Ent{u_t^{(N)}}{u_t^{\otimes N}}}{N} \leq e^{-\frac{1}{2\CLSI}t}  \frac{\Ent{u_0^{(N)}}{u_0^{\otimes N}}}{N} + \frac{\Constinit}{N}\left( 1 + e^{-\gamma t} N^{1-\frac{1}{p}}\mathbf{W}_{2p}^2\left(u_0^{(N)},u_0^{\otimes N}\right)\right),
    \end{equation*}
choosing $p = 2$ yields the fastest rate as the quantity $\mathbf{W}_{2p}^2N^{-\frac{1}{p}}$ is increasing in $p$.
\end{eproof}

\textbf{Acknowledgments.} I warmly thank Antoine Diez and Louis-Pierre Chaintron for their valuable comments and helpful discussions during the preparation of this work. I also would like to express my gratitude to my supervisors, François Bolley and Paul-Éric Chaudru de Raynal, for their guidance, advice and careful readings.

\noindent

\newcommand{\etalchar}[1]{$^{#1}$}


\end{document}